\documentclass{scrartcl}
\usepackage[numbers]{natbib}
\usepackage{fontenc}
\usepackage{shadethm}
\usepackage{amsthm}
\usepackage{float}
\usepackage{bm}
\usepackage{stmaryrd}
\usepackage{mathrsfs} 
\usepackage{amsmath}
\usepackage{amssymb}
\usepackage{wasysym}
\usepackage{sectsty}
\usepackage{hyperref}
\usepackage[a4paper, left=2.7cm, right=2.7cm, top=2.5cm]{geometry}
\usepackage{chngcntr}
\counterwithin*{equation}{section}

\usepackage{cleveref}
\DeclareMathAlphabet{\mathpzc}{OT1}{pzc}{m}{it}

\sectionfont{\normalsize\centering}
\subsectionfont{\footnotesize\centering}
\newcommand{\norm}[1]{\left\lVert#1\right\rVert}

\theoremstyle{plain}
\newtheorem{thm}{Theorem}[section] 

\theoremstyle{definition}
\newtheorem{exmp}[thm]{Example} 
\newtheorem{lem}[thm]{Lemma}
\newtheorem{prop}[thm]{Proposition}
\newtheorem{rem}[thm]{Remark}
\newtheorem{cor}[thm]{Corollary}

\def\XXint#1#2#3{{\setbox0=\hbox{$#1{#2#3}{\int}$ }
		\vcenter{\hbox{$#2#3$ }}\kern-.6\wd0}}

\usepackage[utf8]{inputenc}
\usepackage[german,english,russian]{babel}

\newcounter{MPequ}
\newenvironment{MPEquation}
{\stepcounter{MPequ}%
	\addtocounter{equation}{0}%
	\equation}
{\endequation}

\newcounter{AppA}
\newenvironment{AppA}
{\stepcounter{AppA}%
	\addtocounter{equation}{0}%
	\equation}
{\endequation}

\begin{document}\selectlanguage{english}
\begin{center}
\normalsize \textbf{\textsf{Existence and structure of symmetric Beltrami flows on compact $3$-manifolds}}
\end{center}
\begin{center}
	Wadim Gerner\footnote{\textit{E-mail address:} \href{mailto:gerner@eddy.rwth-aachen.de}{gerner@eddy.rwth-aachen.de}}
\end{center}
\begin{center}
{\footnotesize	RWTH Aachen University, Lehrstuhl f\"ur Angewandte Analysis, Turmstra{\ss}e 46, D-52064 Aachen, Germany}
\end{center}
{\small \textbf{Abstract:} We show that for almost every given symmetry transformation of a Riemannian manifold there exists an eigenvector field of the curl operator, corresponding to a non-zero eigenvalue, which obeys the symmetry. More precisely, given a smooth, compact, oriented Riemannian $3$-manifold $(\bar{M},g)$ with (possibly empty) boundary and a smooth flow of isometries $\phi_t:\bar{M}\rightarrow \bar{M}$ we show that, if $\bar{M}$ has non-empty boundary or if the  infinitesimal generator is not purely harmonic, there is a smooth vector field $\bm{X}$, tangent to the boundary, which is an eigenfield of curl and satisfies $(\phi_t)_{*}\bm{X}=\bm{X}$, i.e. is invariant under the pushforward of the symmetry transformation.
We then proceed to show that if the quantities involved are real analytic and $(\bar{M},g)$ has non-empty boundary, then Arnold's structure theorem applies to all eigenfields of curl, which obey a symmetry and appropriate boundary conditions. More generally we show that the structure theorem applies to all real analytic vector fields of non-vanishing helicity which obey some nontrivial symmetry. A byproduct of our proof is a characterisation of the flows of real analytic Killing fields on compact, connected, orientable $3$-manifolds with and without boundary.
\section{Introduction}
Eigenvector fields of the curl operator corresponding to non-zero eigenvalues, often referred to as (strong) Beltrami fields, are of particular interest in ideal magnetohydrodynamics as well as hydrodynamics, \cite[Chapter II \& Chapter III]{AK98}, where these types of vector fields appear as solutions of the equations of ideal magnetohydrodynamics for the case of constant pressure and a resting fluid, and as solutions of the incompressible Euler equations for appropriate pressure functions. The stationary, incompressible Euler equations may be written as
\begin{equation}
\label{I1}
\bm{v}\times \text{curl}(\bm{v})=\nabla f\quad \text{ and }\quad \text{div}(\bm{v})=0,
\end{equation} 
where $\bm{v}$ is the velocity vector field and $f$ is the Bernoulli function of the system. When it comes to understanding the dynamics of solutions of steady Euler flows, i.e. solutions $\bm{v}$ of (\ref{I1}) for some given function $f$, then Arnold's structure theorem, \cite{A66}, \cite{A74}, \cite[Chapter II, Theorem 1.2]{AK98}, \cite{GK94} asserts that if the quantities involved are all real analytic and $\bm{v}$ and its curl are linearly independent in at least one point, then after removing an appropriate semianalytic subset of codimension at least $1$, the remaining set decomposes into finitely many connected components in each of which the flow of $\bm{v}$ shows a standard behaviour. Arnold's theorem provides a characterisation of steady Euler flows which are not everywhere collinear with their curl. But strong Beltrami fields fail to fall under this standard characterisation and may have very complicated flows, see for instance \cite{DFGHMS86} where the flow of the so called ABC flows is studied in detail. In fact the presence of a 'chaotic' field line in a non-singular, steady Euler flow implies that this flow is Beltrami, \cite[Chapter II, Proposition 6.2]{AK98}. As for the question of how topologically complicated the field line behaviour of Beltrami fields can be see further \cite{EtGr00b},\cite{EP12},\cite{EP15}, where the authors for example show that there exists a Beltrami field on $\mathbb{R}^3$ which admits field line configurations of all (tame) knot and (locally finite) link types simultaneously. In the present work we pose the following question:
\newline
\newline
\textit{Does there exist a class of Beltrami fields, which has a 'well-behaved' flow, i.e. a class of Beltrami fields to which the characterisation of Arnold's structure theorem applies?}
\newline
\newline
It was for instance observed by Cantarella \cite{C99}  that rotationally symmetric Beltrami fields on rotationally symmetric domains diffeomorphic to the solid torus in $\mathbb{R}^3$ show the same behaviour as Arnold's structure theorem describes. Also in \cite{CDGT00} the authors compute all eigenvector fields of curl on the closed unit ball which are tangent to its boundary. It turns out that the Beltrami fields corresponding to the smallest positive eigenvalue are in fact rotationally symmetric and also show the same flow behaviour as steady Euler flows to which Arnold's structure theorem applies. Cantarella seems to have been the first, to the best of my knowledge, who generalised the existence and structure result to symmetric Beltrami flows (even though only for the special case of rotationally symmetric domains).
\newline
One key idea is the well-known fact that the existence of Beltrami fields on abstract manifolds has a variational formulation by means of a helicity constraint $L^2$-energy minimisation \cite{A74},\cite{AL91}, which dates back to Woltjer's work \cite{W58} in 1958. The notion of helicity, as a conserved physical quantity, and its relation to the topology of field lines has been widely studied in the literature, see for example \cite{M69}, \cite{V03}, \cite{MR92}, \cite{CDG00}, \cite{CP10}, \cite{EPT16}. Most notably for us here is the approach by Arnold, \cite{A74}, who used a spectral theoretical argument to prove the existence of global energy minimisers of the constraint minimisation problem, which turn out to be Beltrami fields.
\newline
In this paper we generalise the results of Cantarella to the setting of abstract manifolds with respect to any (continuous) symmetry transformation. Let us point out that Arnold's method allows us to obtain a candidate of a Beltrami field, which may (or a priori may not) obey the symmetry condition. In \cite{C99} the method used in order to show that this candidate is indeed rotationally symmetric makes use of the explicit structure of the Biot-Savart potential of a divergence-free vector field, tangent to the boundary of a domain, and hence is specifically tailored for the Euclidean problem. We take a different path here, which is applicable to all abstract manifolds. First, in order to show the existence of a suitable candidate, we make use of Arnold's approach and the Hodge decomposition for abstract manifolds \cite{S95}, see also \cite{CDG02} for decomposition theorems in $\mathbb{R}^3$. Second, in order to show that this candidate indeed obeys the symmetry condition, we first establish its regularity and then use an approximation argument in combination with an equivalent reformulation of the symmetry property by means of an integral equation, exploiting our specific choice of boundary conditions.
\newline
As for the structure of rotationally symmetric Beltrami flows, Cantarella considers the vector field in cylindrical coordinates and reduces the property of being an eigenvector field of curl to the problem of solving a suitable scalar elliptic Dirichlet problem on appropriate cross sections of the torus, \cite[Proposition 3, 4 and 5]{C99}. In contrast to that we will not study reduced elliptic problems for that matter, but instead follow the exposition of Arnold's proof of the structure theorem. In order to do so we will see that the extra symmetry of a Beltrami field $\bm{X}$ allows us to express the cross product $\bm{X}\times \bm{Y}$, where $\bm{Y}$ is the Killing vector field generating the symmetry, as a gradient field and that since both vector fields commute, $[\bm{X},\bm{Y}]=0$, Arnold's reasoning carries over to our specific situation. A byproduct of our approach is the fact that the structure theorem applies to all non-trivial, real analytic Killing vector fields, as long as there exists a corresponding Beltrami field obeying the induced symmetry. This in combination with our existence result in particular shows that the structure theorem applies to every non-trivial, real analytic Killing field on a compact, oriented, connected Riemannian manifold with non-empty boundary.
\newline
Note that Beltrami fields with non-constant proportionality function are also considered in the literature, so called weak Beltrami fields, see for example \cite{EP16}, \cite{N14}. Sometimes the specification weak and strong is dropped and one has to pay attention which kind of Beltrami fields are under consideration. In the present paper we will be always referring to strong Beltrami fields.
\section{Main results}
\textbf{Conventions:} All manifolds in question are assumed to be Hausdorff, second countable, oriented, smooth, with or without boundary. We will simply write: 'Let $(\bar{M},g)$ be a $3$-manifold' to indicate that $\bar{M}$ has all these properties, is $3$-dimensional and that $g$ is a smooth metric on $\bar{M}$. We call a smooth manifold $\bar{M}$ with or without boundary analytic, if it is additionally equipped with a real analytic atlas for the interior which is compatible with the induced smooth structure. A metric $g$ on a real analytic manifold $\bar{M}$ is said to be real analytic if it is smooth up to the boundary and its restriction to the interior is real analytic with respect to the fixed real analytic structure. We will say 'Let $(\bar{M},g)$ be a real analytic $3$-manifold' and mean that it is a $3$-manifold in the previously defined meaning, that $\bar{M}$ is equipped with a real analytic structure of the interior and that $g$ is real analytic with respect to this structure. Similarly we call a vector field $\bm{X}$ on a real analytic $3$-manifold real analytic if it is smooth up to the boundary and real analytic on the interior.
\newline
We denote by $\mathcal{V}(\bar{M})$ the set of all smooth vector fields on $\bar{M}$, by $\mathcal{V}_n(\bar{M})$ the set of all smooth vector fields $\bm{X}$ such that there exists some $\bm{A}\in \mathcal{V}(\bar{M})$ with $\bm{A}\perp \partial\bar{M}$ and $\text{curl}(\bm{A})=\bm{X}$ and by $\mathcal{K}(\bar{M})$ the set of smooth Killing fields on $\bar{M}$ generating global isometries, i.e. $\bm{Y}\in \mathcal{K}(\bar{M})$ if and only if $\bm{Y} \parallel \partial\bar{M}$, $\bm{Y}\in \mathcal{V}(\bar{M})$ and $\bm{Y}$ gives rise to a global flow $\phi_t$ such that $\phi_t:(\bar{M},g)\rightarrow (\bar{M},g)$ are all isometries. In our applications the manifolds will be compact so that flows will be automatically global. Similarly we denote by $\mathcal{V}^{\omega}(\bar{M})$ and $\mathcal{K}^{\omega}(\bar{M})$ the corresponding real analytic spaces if $(\bar{M},g)$ is real analytic. Note that smooth Killing fields on real analytic manifolds are necessarily real analytic. Nonetheless we will use distinct notations to emphasise the setting we are currently working in.
\newline
A (strong) Beltrami field is a smooth vector field $\bm{X}\in \mathcal{V}(\bar{M})$ which is an eigenvector field of the curl operator corresponding to a non-zero eigenvalue (in particular $\bm{X}\not\equiv 0$).
\newline
Lastly, given a smooth $3$-manifold $(\bar{M},g)$ and any fixed $\bm{Y}\in \mathcal{K}(\bar{M})$, we define $\mathcal{V}^{\bm{Y}}_n(\bar{M}):=\{\bm{X}\in \mathcal{V}_n(\bar{M})|[\bm{X},\bm{Y}]\equiv 0 \}$, where $[\cdot,\cdot]$ denotes the Lie bracket of vector fields.
\begin{thm}[Conditional existence]
\label{MT1}
Let $(\bar{M},g)$ be a compact $3$-manifold and $\bm{Y}\in \mathcal{K}(\bar{M})$. If $\mathcal{V}^{\bm{Y}}_n(\bar{M})\neq \{0\}$, then there exists a vector field $\bm{X}\in \mathcal{V}_n(\bar{M})\setminus \{0\}$ and some $\lambda\in \mathbb{R}\setminus \{0\}$ such that
\begin{equation}
\label{M1}
\text{curl}(\bm{X})=\lambda \bm{X}\text{ and }[\bm{Y},\bm{X}]\equiv 0.
\end{equation}
\end{thm}
\begin{rem}
\label{MR2}
\begin{enumerate}
\item Every smooth flow of global isometries is generated by a Killing field and the condition $[\bm{Y},\bm{X}]\equiv 0$ is obviously equivalent to the statement $(\phi_t)_{*}\bm{X}=\bm{X}$ for all times $t$, where $\phi_t$ denotes the flow of $\bm{Y}$.
\item It follows from elliptic estimates that if $(\bar{M},g)$ is a real analytic $3$-manifold, then $\bm{X}$ is real analytic, which is of importance in view of the upcoming structure theorems.
\item Our assumption on the Killing field $\bm{Y}$ to be tangent to the boundary is crucial for the existence proof. In particular our arguments do not generalise to the case of vector fields $\bm{Y}$ which only satisfy the Killing equations, but are not necessarily tangent to the boundary.
\item We will see that we have the relation $\text{grad}(g(\bm{X},\bm{Y}))=\lambda \bm{Y}\times \bm{X}$ for every Beltrami field $\bm{X}$ which commutes with a given Killing field $\bm{Y}\in \mathcal{K}(\bar{M})$. Thus, if $\bm{X}$ is tangent to the boundary, the tangent part of $\text{grad}(g(\bm{X},\bm{Y}))$ vanishes and hence $g(\bm{X},\bm{Y})|_{\partial\bar{M}}$ is a locally constant function. In particular on any connected component of the boundary, which is not diffeomorphic to the torus, the vector field $\bm{Y}$, being tangent to the boundary, must vanish in at least one point and thus on each such component $\bm{X}$ and $\bm{Y}$ are everywhere $g$-orthogonal. The flow of Killing fields in dimension $2$ is well understood, see also \cref{Appendix} of the present paper, which gives us a qualitative understanding of the boundary field line behaviour of symmetric Beltrami fields, since the restriction of Killing fields to the boundary are again Killing fields with respect to the pullback metric. In particular if $\text{grad}_{B}(g(\bm{Y},\bm{Y}))$ denotes the gradient of the function $g|_{\partial\bar{M}}(\bm{Y},\bm{Y})$ with respect to the pullback metric on the boundary, then $\text{grad}_{B}(g(\bm{Y},\bm{Y}))$ and $\bm{X}$ are everywhere linearly dependent on the considered boundary component.
\end{enumerate}
\end{rem}
The following result presents a broad class of situations in which the condition $\mathcal{V}^{\bm{Y}}_n(\bar{M})\neq \{0\}$ is satisfied.
\begin{prop}
\label{ExtraP3}
Let $(\bar{M},g)$ be a compact $3$-manifold and $\bm{Y}\in \mathcal{K}(\bar{M})$. Then $\mathcal{V}^{\bm{Y}}_n(\bar{M})\neq \{0\}$ whenever one of the following conditions is satisfied
\begin{enumerate}
\item $\partial\bar{M}\neq \emptyset$.
\item There exists a non-harmonic Killing field on $\bar{M}$.
\item $\bm{Y}\equiv 0$.
\end{enumerate}
\end{prop}
The second condition more precisely demands that there exists some $\tilde{\bm{Y}}\in \mathcal{K}(\bar{M})$ (which may, but need not to coincide with the Killing field $\bm{Y}$ of interest) which satisfies $\text{curl}(\tilde{\bm{Y}})\not\equiv 0$. The remaining case is dealt with in the following proposition
\begin{prop}
\label{Extrap4}
Let $(M,g)$ be a compact, connected $3$-manifold with empty boundary and suppose $\bm{Y}\in \mathcal{K}(M)\setminus\{0\}$ satisfies $\text{curl}(\bm{Y})\equiv 0$. Then the following four statements are equivalent
\begin{enumerate}
\item There exists some $f\in C^{\infty}(M)$ and $\lambda >0$ such that $\bm{X}:=\bm{Y}\times \text{grad}(f)-f\lambda\bm{Y}$ satisfies
\[
\text{curl}(\bm{X})=\lambda\bm{X},\text{ }[\bm{X},\bm{Y}]\equiv 0,\text{ and }\bm{Y}\times \bm{X}\not\equiv 0
\]
\item There is some $\bm{X}\in \mathcal{V}_n(M)\setminus\{0\}$ and $\lambda \in \mathbb{R}\setminus \{0\}$ with $\text{curl}(\bm{X})=\lambda \bm{X}\text{ and }[\bm{X},\bm{Y}]\equiv 0$
\item $\mathcal{V}_n^{\bm{Y}}(M)\neq \{0\}$
\item $\{f\in C^{\infty}(M)|df\not\equiv 0\text{ and }df(\bm{Y})\equiv 0\}\neq \emptyset$
\end{enumerate}
\end{prop}
Note that in the last statement we have the identities $df(\bm{Y})=g(\bm{Y},\text{grad}(f))=\bm{Y}(f)=\mathcal{L}_{\bm{Y}}(f)$, where the latter denotes the Lie derivative along $\bm{Y}$. Thus there exists a symmetry obeying Beltrami field if and only if the considered Killing field admits a nontrivial first integral. The following example shows that not every Killing field admits a nontrivial first integral.
\begin{exmp}
\label{ExtraE5}
Consider the $3$-torus with its flat metric $(T^3,g_{\text{F}})$. Viewing $T^3$ as a unit cube in $\mathbb{R}^3$ with opposite faces identified, the unit vector fields, inducing translations along the coordinate axes on $\mathbb{R}^3$, descend to well-defined vector fields on $T^3$, which we denote by $e_1,e_2,e_3$ respectively. It is not hard to check that $\mathcal{K}(T^3,g_{\text{F}})$ is spanned by these $3$ vector fields and that all of these vector fields are irrotational, i.e. $\text{curl}(e_i)\equiv 0$, so that the results of \cref{ExtraP3} do not apply. Then $\bm{Y}:=e_1+\sqrt{2}e_2+\sqrt{6}e_3$ is also a Killing field as a linear combination of such vector fields and further $(1,\sqrt{2},\sqrt{6})$ are rationally independent, i.e. $k_1+\sqrt{2}k_2+\sqrt{6}k_3=0$ for $k_1,k_2,k_3\in \mathbb{Z}$ implies that $k_1=k_2=k_3=0$. It follows that the images of all field lines of $\bm{Y}$ are dense in $T^3$, \cite[Proposition 1.5.1]{HasKa95}. Note that the condition $\mathcal{L}_{\bm{Y}}(f)\equiv 0$ is equivalent to the statement that $f$ is constant along every field line of $\bm{Y}$. Thus, since every field line of $\bm{Y}$ is dense in $T^3$, we see that any smooth function $f$ which satisfies $\mathcal{L}_{\bm{Y}}(f)\equiv 0$ is constant on some dense subset of $T^3$ and hence is constant on all of $T^3$, i.e. $\{f\in C^{\infty}(T^3)|df\not\equiv 0\text{ and }df(\bm{Y})\equiv 0\}=\emptyset$, see also \cite[Corollary 1.4.4]{HasKa95}. In particular \cref{Extrap4} shows that $\bm{Y}$ does not admit any symmetry obeying eigenfield of curl, corresponding to a non-zero eigenvalue. Of course $\bm{Y}$ itself is an eigenfield of curl, corresponding to the eigenvalue zero, and commutes with itself.
\end{exmp}
In fact the arguments we provide imply a version with more than one symmetry, as long as the symmetries are compatible.
\begin{cor}
\label{MCA1}
Let $(\bar{M},g)$ be a compact $3$-manifold. Suppose $N\in \mathbb{N}$ and $\bm{Y}_1,\dots,\bm{Y}_N\in \mathcal{K}(\bar{M})$ are Killing fields such that $[\bm{Y}_1,\bm{Y}_j]\equiv 0$ for all $1\leq j\leq N$. If $\partial\bar{M}\neq \emptyset$ or if $\text{curl}(\bm{Y}_1)$ is not identical zero, then there exists a vector field $\bm{X}\in \mathcal{V}_n(\bar{M})\setminus \{0\}$ and some $\lambda \in \mathbb{R}\setminus\{0\}$ with
\begin{equation}
\text{curl}(\bm{X})=\lambda \bm{X}\text{ and }[\bm{X},\bm{Y}_j]\equiv 0\text{ for all }1\leq j\leq N.
\end{equation}
\end{cor}
Next we want to state a structure theorem for symmetric Beltrami fields. Let us give an informal description of its implication first to get an intuition for the upcoming result. The structure theorem tells us roughly that it is possible (under mild assumptions) to extract for every fixed Beltrami field, which obeys a symmetry, an at most two dimensional subset of the manifold, such that the complement of this subset consists of three dimensional chambers, each of which has a standard behaviour with regards to the flow of the considered Beltrami field. Namely, each such chamber is diffeomorphic to $\mathbb{R}\times T^2$, i.e. consists of layers of tori, all of which are invariant under the flow of the Beltrami field, i.e. each field line of the Beltrami field starting at such a torus is confined to the torus for all times. Moreover on every such fixed torus all field lines on the torus show the same behaviour, namely either all field lines are closed or all field lines are dense in the torus, which provides us with a qualitative understanding of the flow on these tori. 
\newline
\newline
We denote by $\mathcal{H}^{\alpha}$ for $0\leq \alpha$ the corresponding Hausdorff measure induced by the standard metric on a given connected Riemannian manifold, \cite[Theorem 2.55]{L18}. We will also use the terminology of Federer, regarding the rectifiability of sets in metric spaces \cite[Definition 3.2.14]{Fed69}
\begin{thm}[Structure theorem]
\label{MT3}
Let $(\bar{M},g)$ be a compact, connected real analytic $3$-manifold, $Y\in \mathcal{K}^{\omega}(\bar{M})$ be a Killing field and $\bm{X}\in \mathcal{V}^{\omega}(\bar{M})$ be a Beltrami field which is tangent to the boundary of $\bar{M}$. If $\bm{X}$ satisfies $[\bm{Y},\bm{X}]\equiv 0$ and $\bm{X}$ and $\bm{Y}$ are in at least one point linearly independent, then there exists a compact, $\mathcal{H}^2$-countably-$2$ rectifiable subset $\Gamma\subset \bar{M}$ in the sense of Federer such that
\begin{enumerate}
\item $\partial\bar{M}\subseteq \Gamma$.
\item $\bar{M}\setminus \Gamma$ is the disjoint union of tori $T_i$ which are real analytically embedded in the interior $M:=\text{int}(\bar{M})$ of $\bar{M}$ and which are invariant under the flow of $\bm{X}$.
\item On each $T_i$ the orbits of $\bm{X}$ starting at $T_i$ are either all non-constant and closed (all of the same period) or the image of every orbit starting at $T_i$ is dense in $T_i$.
\item For every $T_i$, the connected component $U_i$ of $\bar{M}\setminus \Gamma$ containing it, is an open neighbourhood of $T_i$ which is invariant under the flow of $\bm{X}$ and there is a diffeomorphism $\psi:(-\epsilon_1,\epsilon_2)\times T_i\rightarrow U_i$ for suitable constants $\epsilon_2,\epsilon_1>0$, which satisfies $\psi|_{\{0\}\times T_i}=\text{Id}|_{T_i}$ and for each fixed $-\epsilon_1<t<\epsilon_2$, $\psi(t,T_i)=T_{j(t)}$ for some suitable invariant torus $T_{j(t)}$. In particular every connected component of $\bar{M}\setminus \Gamma$ is diffeomorphic to $\mathbb{R}\times T^2$, where $T^2$ denotes the $2$-torus.
\item $\Gamma$ is invariant under the flow of $\bm{X}$.
\end{enumerate}
\end{thm}
\begin{rem}
\label{MR4}
\begin{enumerate}
\item  The set $\Gamma$ has a Hausdorff dimension of at most $2$, hence is a null set and nowhere dense. In fact $\Gamma$ is the union of $\partial\bar{M}$ and at most countably many real analytically embedded (into the interior $M$ of $\bar{M}$) surfaces, curves and points.
\item If $\partial\bar{M}=\emptyset$, then $\Gamma$ is semianalytic and $\bar{M}\setminus \Gamma$ has at most finitely many connected components, see also \cite{BM88} for an introduction into this topic.
\item Even if we additionally demand $\bm{Y}\not\equiv 0$ in \cref{MT3}, then $\bm{X}$ and $\bm{Y}$ can turn out to be everywhere collinear. To see this consider the Hopf vector field $\bm{H}$ on $S^3\subset \mathbb{R}^4$ (equipped with the round metric) given by $(x_1,x_2,x_3,x_4)\mapsto (-x_2,x_1,-x_4,x_3)$. This vector field is both, Beltrami and Killing, as can be easily checked. Thus the Hopf field is a Beltrami field and commutes with itself, hence admits a symmetry, but satisfies $\bm{H}\times \bm{H}\equiv 0$. Therefore we cannot simply set $\bm{X}=\bm{H}=\bm{Y}$ in \cref{MT3} to conclude that the structure theorem applies. However, it can be verified that the structure theorem applies nonetheless in this case by explicitly computing the flow of $\bm{H}$. The deeper reason as to why the structure theorem applies anyway, lies in the fact that the Hopf field admits a nontrivial first integral. For instance, the map $f:S^3\rightarrow \mathbb{R}, (x_1,x_2,x_3,x_4)\mapsto x^2_1+x^2_2$ is invariant under the flow of the Hopf field, so that \cref{MTE2} applies.
\item Since $\bm{Y}\in \mathcal{K}^{\omega}(\bar{M})$ is tangent to the boundary, it generates a global flow, and we can apply the exact same reasoning of the proof of \cref{MT3} to the flow of $\bm{Y}$, i.e. if the requirements of \cref{MT3} are satisfied then for the same set $\Gamma$ as in the case of $\bm{X}$, the tori of $\bar{M}\setminus \Gamma$ which are invariant under $\bm{X}$ are also invariant under the flow of $\bm{Y}$.
\item Given a vector field $\bm{X}\in \mathcal{V}_n(\bar{M})$ we define its helicity $\mathcal{H}(\bm{X})$ as the $L^2$-inner product of $\bm{X}$ and any of its vector potentials, where the value turns out to be independent of the choice of potential. If $\bm{X}\in \mathcal{V}_n(\bar{M})$ is real analytic with a non-zero helicity and commutes with some real analytic Killing field $\bm{Y}\in \mathcal{K}^{\omega}(\bar{M})\setminus \{0\}$, then the structure theorem applies to $\bm{X}$, whenever $\bar{M}$ has non-empty boundary. Thus more generally every symmetric, real analytic vector field with non-vanishing helicity on a compact, connected $3$-manifold with non-empty boundary obeys the structure theorem, see also \cref{PRP1}. Note however that \cref{MT3} applies to all Beltrami fields, tangent to the boundary, even those which are not necessarily contained in $\mathcal{V}_n(\bar{M})$.
\item In Arnold's classical result \cite[Chapter II, Theorem 1.2]{AK98} regarding steady Euler flows in the presence of boundary, the chambers can be of two types. Either they consist of layers of tori, which are invariant under the flow or they fibre into layers of annuli, diffeomorphic to $\mathbb{R}\times S^1$, which are invariant under the flow. In contrast to that, in our situation, the chambers always consist of layers of tori. This is due to the fact that in our case the involved vector fields, the Beltrami as well as Killing field, are tangent to the boundary, whereas in the classical situation, the velocity vector field $\bm{v}$ is tangent to the boundary, while its curl must not be, which implies that there are more situations to consider in the classical structure theorem for Euler flows.
\end{enumerate}
\end{rem}
We already gave the example of the Hopf vector field which is simultaneously Beltrami and Killing. However in the case of non-empty boundary this cannot happen
\begin{lem}
\label{ML5}
Let $(\bar{M},g)$ be a compact, connected $3$-manifold with non-empty boundary. Let $\bm{X}\in \mathcal{V}_n(\bar{M})$ be a Beltrami field and $\bm{Y}\in \mathcal{K}(\bar{M})$. If $\bm{Y}\times \bm{X}\equiv 0$, then $\bm{Y}\equiv 0$.
\end{lem}
\begin{rem}
\label{MR6}
If $\bm{Y}\in \mathcal{K}^{\omega}(\bar{M})$ with  $\bm{Y}\not\equiv 0$ and $(\bar{M},g)$ is real analytic and with non-empty boundary, then \cref{MT1}, \cref{ExtraP3} and \cref{MR2} imply that the corresponding Beltrami field $\bm{X}$ whose existence is established is real analytic (and obviously also tangent to the boundary), while \cref{ML5} implies that the requirements of \cref{MT3} are satisfied and hence in this case the conclusions of the structure theorem always apply to the Beltrami field $\bm{X}$.
\end{rem}
Note also that according to \cref{MR4} the structure theorem applies to $\bm{Y}$ whenever it applies to the corresponding Beltrami field. Thus it follows from \cref{MR6} that we have in particular the following, see also \cref{MP7}:
\begin{thm}[Characterisation of real analytic Killing flows, non-empty boundary case]
\label{MCE1}
Let $(\bar{M},g)$ be a compact, connected real analytic $3$-manifold with non-empty boundary and let $\bm{Y}\in \mathcal{K}^{\omega}(\bar{M})$ be a real analytic Killing field. If $\bm{Y}\not\equiv 0$, then the conclusions of the structure \cref{MT3} apply to $\bm{Y}$. If further $\bm{Y}\times\text{curl}(\bm{Y})\equiv 0$, then after rescaling $\bm{Y}$ by a suitable constant factor $c>0$, all field lines of the rescaled version of $\bm{Y}$ are unit speed geodesics.
\end{thm}
The boundaryless case is more diverse because we do not have \cref{ML5} at hand and have to deal with more possible situations. A proof based on the existence of symmetric Beltrami flows is given in a later section.
\begin{thm}[Characterisation of real analytic Killing flows, empty boundary case]
\label{MTE2}
Let $(M,g)$ be a compact, connected, real analytic $3$-manifold with empty boundary and let $\bm{Y}\in \mathcal{K}^{\omega}(M)$ be a real analytic Killing field. Then exactly one of the following $2$ situations occurs
\begin{enumerate}
\item There does not exist any $\kappa \in \mathbb{R}$ with $\text{curl}(\bm{Y})=\kappa \bm{Y}$: In this case the conclusions of the structure \cref{MT3} apply to $\bm{Y}$.
\item There exists some $\kappa\in \mathbb{R}$ with $\text{curl}(\bm{Y})=\kappa\bm{Y}$: In this case either $\bm{Y}\equiv 0$, or otherwise, after rescaling $\bm{Y}$ by a suitable constant factor $c>0$, all field lines of the rescaled version of $\bm{Y}$ are unit speed geodesics. If in the latter case, the Killing field $\bm{Y}$ admits a nontrivial first integral, i.e. if there exists a non-constant $C^1$-function $f:M\rightarrow \mathbb{R}$, such that $\bm{Y}(f)\equiv 0$, then the conclusions of the structure \cref{MT3} apply to $\bm{Y}$.
\end{enumerate}
Moreover in all these cases the sets $\Gamma$ as well as $M\setminus \Gamma$ obtained from the structure theorem are semianalytic and both have finitely many connected components.
\end{thm}
Note that it is enough to demand $C^1$-regularity of the first integral. The proof we provide guarantees under these assumptions the existence of a nontrivial, real analytic first integral. Observe also that the Killing field $\bm{Y}$, which we defined in \cref{ExtraE5}, is not identically zero and provides an example for a Killing field to which the structure theorem does not apply. Indeed if the structure theorem were to apply, then the set $T^3\setminus \Gamma$ is non-empty and is the disjoint union of tori, each of which is invariant under the flow of $\bm{Y}$. Thus, in this case, $\bm{Y}$ must in particular admit at least one field line which is confined to an embedded $2$-torus and hence its image cannot be dense in $T^3$, which is a contradiction.
\newline
The following proposition shows that curl-free Killing fields as well as Killing fields which are Beltrami fields are always of a very special kind
\begin{prop}
\label{MP7}
Let $(\bar{M},g)$ be a $3$-manifold and $\bm{Y}\in \mathcal{K}(\bar{M})$. Suppose that either $\text{curl}(\bm{Y})\equiv 0$ or that $\bm{Y}$ is a Beltrami field, then
\begin{enumerate}
\item $\text{grad}(g(\bm{Y},\bm{Y}))\equiv 0$
\item $\nabla_{\bm{Y}}\bm{Y}\equiv 0$, i.e. all field lines of $\bm{Y}$ are constant speed geodesics.
\end{enumerate}
If in addition $\partial\bar{M}\neq \emptyset$ and it has at least one compact connected component which is not homeomorphic to the torus and $\bar{M}$ is connected, then any $\bm{Y}\in \mathcal{K}(\bar{M})$ which is everywhere collinear with some Beltrami field or is irrotational satisfies $\bm{Y}\equiv 0$. 
\end{prop}
The last part of \cref{MP7} in combination with the fact that boundaries of compact manifolds are compact is of particular interest in view of the requirements of the structure theorem. Note that this conclusion applies to any Beltrami field in $\mathcal{V}(\bar{M})$ and thus in particular to Beltrami fields which are tangent to the boundary but possibly not in $\mathcal{V}_n(\bar{M})$, in which case \cref{ML5} is not at hand.
\begin{exmp}
\label{ME8}
Let $\bm{Y}:\mathbb{R}^3\rightarrow \mathbb{R}^3, (x,y,z)\mapsto (-y,x,0)$, then $\bm{Y}$ is a real analytic vector field whose flow $\phi_t$ induces rotations around the $z$-axis and hence is a family of isometries (with respect to the Euclidean metric). If we let $\Omega \subset \mathbb{R}^3$ be any nonempty domain with smooth boundary, which is invariant under the flow, i.e. $\phi_t(\Omega)=\Omega$, then the restriction of $\bm{Y}$ to $\Omega$ defines a Killing field, tangent to $\partial\Omega$. Since $\Omega$ has non-empty boundary, \cref{MT1} and \cref{ExtraP3} show the existence of rotationally symmetric Beltrami fields on each such domain. On the other hand if $\bm{X}$ is any rotationally symmetric (around the $z$-axis) real analytic Beltrami field, then $\bm{Y}$ and $\bm{X}$ cannot be everywhere collinear. Assume they are, then since $\bm{X}$ is not the zero vector field the set of points $U$ at which $\bm{X}=(X^1,X^2,X^3)$ does not vanish in $\Omega$ is a non-empty open subset of $\Omega$. Thus on $U$ there is a smooth function $f:U\rightarrow \mathbb{R}$ such that $\bm{Y}=f\bm{X}$, which immediately implies $X^3= 0$ on $U$ whenever $f$ does not vanish. Further we conclude $(0,0,2)=\text{curl}(\bm{Y})=\nabla f\times \bm{X}+f\lambda \bm{X}$ where $\lambda\neq 0$ is the corresponding eigenvalue. Multiply this equation by $\bm{X}$ at any potential point at which $f$ does not vanish to obtain $0=2X^3=(0,0,2)\cdot \bm{X}=f\lambda |\bm{X}|_2^2$, with $\cdot$ denoting the standard inner product. We have $\lambda \neq 0$ and $|\bm{X}|_2\neq 0$ on $U$ by definition of $U$. Thus $f\equiv 0$ on $U$, which in turn implies that $\bm{Y}$ vanishes on some open set, which is a contradiction because $\bm{Y}$ vanishes exactly along the $z$-axis. Thus the structure \cref{MT3} applies to every rotationally symmetric Beltrami field, tangent to the boundary (not only to elements of $\mathcal{V}_n(\bar{\Omega})$ and even if $\Omega$ is for instance a rotationally symmetric solid torus, i.e. even if \cref{MP7} does not apply). This recovers the results from \cite{C99}. Note however that Cantarella studies the flow structure in more detail for the case of rotationally symmetric solid tori. For instance he shows that the rotationally symmetric eigenfield corresponding to the smallest positive eigenvalue does not vanish anywhere on $\Omega$ or its boundary \cite[Theorem 7]{C99}.
\end{exmp}
\section{Proof of \cref{MT1}}
The following elementary fact is the cornerstone of the existence and structure theorems
\begin{lem}
\label{PT1L1}
Let $(\bar{M},g)$ be a $3$-manifold and $\bm{Y}\in \mathcal{K}(\bar{M})$, then $\text{div}(\bm{Y})\equiv 0$ and for every $\bm{X}\in \mathcal{V}(\bar{M})$ we have the identity
\begin{equation}
\label{PT11}
\text{grad}(g(\bm{X},\bm{Y}))=\bm{Y}\times \text{curl}(\bm{X})+[\bm{Y},\bm{X}].
\end{equation}
\end{lem}
\underline{Proof of \cref{PT1L1}:} The proof follows from direct calculations and the Killing equations. Working in normal coordinates around interior points simplifies the calculations for the interior and a density argument extends it to all of the manifold. $\square$
\newline
\newline
We make the following definitions 
\begin{equation*}
\mathcal{H}(\bar{M}):=\{\bm{\Gamma}\in \mathcal{V}(\bar{M})| \text{div}(\bm{\Gamma})\equiv 0\text{ and } \text{curl}(\bm{\Gamma})\equiv 0 \}\text{, }\mathcal{H}_N(\bar{M}):=\{\bm{\Gamma}\in \mathcal{H}(\bar{M})| \bm{\Gamma}\parallel \partial\bar{M} \}
\end{equation*}
\begin{equation}
\label{PT12}
\text{ and }\mathcal{H}_D(\bar{M}):=\{\bm{\Gamma}\in \mathcal{H}(\bar{M})| \bm{\Gamma}\perp \partial\bar{M} \}
\end{equation}
\begin{lem}
\label{PT1L2}
Let $(\bar{M},g)$ be a compact $3$-manifold. Then for every $\bm{Y}\in \mathcal{K}(\bar{M})$ and every $\bm{\Gamma}\in \mathcal{H}_N(\bar{M})\cup \mathcal{H}_D(\bar{M})$ we have $[\bm{Y},\bm{\Gamma}]\equiv 0$. 
\end{lem}
\underline{Proof of \cref{PT1L2}:} According to \cref{PT1L1} and since $\bm{\Gamma}$ is irrotational we have  $\text{grad}(g(\bm{Y},\bm{\Gamma}))=[\bm{Y},\bm{\Gamma}]$. We conclude $\text{curl}([\bm{Y},\bm{\Gamma}])\equiv 0$. On the other hand since $\bm{\Gamma}$ as well as $\bm{Y}$ are both divergence-free we have by standard vector calculus identities $[\bm{Y},\bm{\Gamma}]=\text{curl}\left(\bm{\Gamma}\times \bm{Y} \right)$ and thus overall $[\bm{Y},\bm{\Gamma}]\in \mathcal{H}(\bar{M})$. Now if $\bm{\Gamma}\in \mathcal{H}_N(\bar{M})$, then $\bm{Y}$ and $\bm{\Gamma}$ are tangent to the boundary which implies that so is their Lie-bracket, i.e. $[\bm{Y},\bm{\Gamma}]\in \mathcal{H}_N(\bar{M})$. But it follows from the integration by parts formula that each gradient field in $\mathcal{H}_N(\bar{M})$ is the zero vector field. On the other hand if $\bm{\Gamma}\in \mathcal{H}_D(\bar{M})$, we see that $g(\bm{Y},\bm{\Gamma})$ vanishes on the boundary, since $\bm{Y}$ is tangent to it, while $\bm{\Gamma}$ is normal to it. We conclude that the tangent part of $\text{grad}(g(\bm{Y},\bm{\Gamma}))=[\bm{Y},\bm{\Gamma}]$ vanishes, i.e. $[\bm{Y},\bm{\Gamma}]\in \mathcal{H}_D(\bar{M})$. But $[\bm{Y},\bm{\Gamma}]$ admits a vector potential, which in combination with the standard integration by parts formula,\cite[Proposition 2.1.2]{S95}, yields $[\bm{Y},\bm{\Gamma}]\equiv 0$. $\square$
\newline
\newline
Recall that $\mathcal{V}_n(\bar{M})$ denotes the set of all smooth vector fields $\bm{X}\in \mathcal{V}(\bar{M})$ which admit a vector potential $\bm{A}\in \mathcal{V}(\bar{M})$, which is normal to the boundary.
\begin{cor}
\label{PT1C3}
Let $(\bar{M},g)$ be a compact $3$-manifold and $\bm{Y}\in \mathcal{K}(\bar{M})$, then there exist unique $\bm{X}\in \mathcal{V}_n(\bar{M})$ and $\bm{\Gamma}\in \mathcal{H}_N(\bar{M})$ such that $\bm{Y}=\bm{X}+\bm{\Gamma}$ and we further have the identity $[\bm{Y},\bm{X}]\equiv 0$.
\end{cor}
\underline{Proof of \cref{PT1C3}:} The decomposition follows immediately from the Hodge-Morrey decomposition \cite[Corollary 3.5.2]{S95} and the fact that $\bm{Y}$ is divergence-free and tangent to the boundary. The remaining claim is a consequence of \cref{PT1L2} and the fact that any vector field commutes with itself. $\square$
\begin{cor}
\label{ExtraT1C4}
Let $(\bar{M},g)$ be a compact, connected $3$-manifold with non-empty boundary. Then every $\bm{Y}\in \mathcal{K}(\bar{M})$ admits a vector potential. More precisely $\bm{Y}=\text{curl}(\bm{B})$ for some suitable $\bm{B}\in \mathcal{V}_n(\bar{M})$.
\end{cor}
\underline{Proof of \cref{ExtraT1C4}:} Since Killing fields are divergence-free, it follows from the Hodge-Morrey-Friedrichs decomposition, \cite[Corollary 3.5.2]{S95}, that we can write $\bm{Y}=\text{curl}(\tilde{\bm{B}})+\bm{\Gamma}$, where $\tilde{\bm{B}}\in \mathcal{V}(\bar{M})$ and $\bm{\Gamma}\in \mathcal{H}_D(\bar{M})$ are suitably chosen smooth vector fields. It then also follows from the Hodge-Morrey decomposition, that the $L^2$-orthogonal projection $\bm{B}$ of $\tilde{\bm{B}}$ onto the space $\mathcal{V}_n(\bar{M})$ has the same curl as $\tilde{\bm{B}}$. Thus $\bm{Y}=\text{curl}(\bm{B})+\bm{\Gamma}$. In view of the fact that the considered decompositions are all $L^2$-orthogonal, it is therefore enough to show that $\mathcal{K}(\bar{M})$ is $L^2$-orthogonal to the space $\mathcal{H}_D(\bar{M})$. To this end suppose $\bm{Y}\in\mathcal{K}(\bar{M})$ and $\bm{\Gamma}\in \mathcal{H}_D(\bar{M})$, then by \cref{PT1L2} we have $[\bm{Y},\bm{\Gamma}]\equiv 0$ and from \cref{PT1L1}, keeping in mind that $\bm{\Gamma}$ is curl-free, we get $[\bm{Y},\bm{\Gamma}]=\text{grad}(g(\bm{Y},\bm{\Gamma}))$. By connectedness of $\bar{M}$, there must be a constant $c\in \mathbb{R}$ with $g(\bm{Y},\bm{\Gamma})\equiv c$. Now $\bm{Y}$ is tangent to the boundary, while $\bm{\Gamma}$ is normal to it. Thus since we assume $\partial\bar{M}\neq \emptyset$ we deduce from this that $c=0$ and consequently $g(\bm{Y},\bm{\Gamma})\equiv 0$. Integrating this equation with respect to the Riemannian volume form gives us the claimed $L^2$-orthogonality property. $\square$
\newline
\newline
Let us now fix some further notation. In the following, $\bm{Y}\in \mathcal{K}(\bar{M})$, is some fixed Killing field
\begin{equation*}
\mathcal{V}_T(\bar{M}):=\{\bm{A}\in \mathcal{V}(\bar{M})| \bm{A}\perp \partial\bar{M}\text{ and }\exists \bm{W}\in \mathcal{V}(\bar{M}): \bm{A}=\text{curl}(\bm{W}) \},
\end{equation*}
\begin{equation}
\label{PT13}
\mathcal{V}^{\bm{Y}}_n(\bar{M}):=\{\bm{X}\in \mathcal{V}_n(\bar{M})|[\bm{Y},\bm{X}]\equiv 0 \}\text{, }\mathcal{V}^{\bm{Y}}_T(\bar{M}):=\{\bm{A}\in \mathcal{V}_T(\bar{M})|[\bm{Y},\bm{A}]\equiv 0 \}.
\end{equation}
By $L^2\mathcal{V}^{\bm{Y}}_n(\bar{M})$ and $H^1\mathcal{V}^{\bm{Y}}_T(\bar{M})$ we denote the corresponding $L^2$ and $H^1$ completions respectively, see \cite[Chapter 1.3]{S95} for an introduction to Sobolev spaces on abstract manifolds.
\begin{lem}
\label{PT1L4}
Let $(\bar{M},g)$ be a compact $3$-manifold and $\bm{Y}\in \mathcal{K}(\bar{M})$, then the following operator is a bounded, linear isomorphism
\begin{equation}
\label{PT14}
\text{curl}: \left(H^1\mathcal{V}^{\bm{Y}}_T(\bar{M}),\norm{\cdot}_{H^1}\right)\rightarrow \left(L^2\mathcal{V}^{\bm{Y}}_n(\bar{M}),\norm{\cdot}_{L^2}\right)
\end{equation}
and the inverse operator $\text{curl}^{-1}$ is $L^2$-compact, i.e. the operator
\begin{equation}
\label{PT15}
\text{curl}^{-1}:\left(L^2\mathcal{V}^{\bm{Y}}_n(\bar{M}),\norm{\cdot}_{L^2}\right)\rightarrow \left(H^1\mathcal{V}^{\bm{Y}}_T(\bar{M}),\norm{\cdot}_{L^2}\right)
\end{equation}
is compact.
\end{lem}
\underline{Proof of \cref{PT1L4}:} First of all we consider the operator $\text{curl}:\mathcal{V}_T(\bar{M})\rightarrow \mathcal{V}_n(\bar{M})$. It is obviously well-defined and injectivity follows immediately, while surjectivity follows from an application of the Hodge-Morrey decomposition. Obviously this operator is $\norm{\cdot}_{H^1}$-$\norm{\cdot}_{L^2}$ continuous and hence has a unique continuous extension to the closures $\text{curl}:\left(H^1\mathcal{V}_T(\bar{M}),\norm{\cdot}_{H^1}\right)\rightarrow \left(L^2\mathcal{V}_n(\bar{M}),\norm{\cdot}_{L^2}\right)$. As for the inverse operator $\text{curl}^{-1}:\mathcal{V}_n(\bar{M})\rightarrow \mathcal{V}_T(\bar{M})$ one can use elliptic estimates, \cite[Lemma 2.4.10]{S95}, to show that this operator is $\norm{\cdot}_{L^2}$-$\norm{\cdot}_{H^1}$ continuous and hence extends to a continuous operator $\text{curl}^{-1}:\left(L^2\mathcal{V}_n(\bar{M}),\norm{\cdot}_{L^2}\right)\rightarrow \left( H^1\mathcal{V}_T(\bar{M}),\norm{\cdot}_{H^1}\right)$. These extended operators are again inverses of one another and the Sobolev embedding theorem \cite[Theorem 1.3.6]{S95} implies the $\norm{\cdot}_{L^2}$ compactness of the extended operator (see also \cite{YG90} where the authors introduce similar boundary conditions on Euclidean domains and show compactness of the corresponding operator). Thus if we can show that $\text{curl}^{-1}\left(L^2\mathcal{V}^{\bm{Y}}_n(\bar{M}) \right)=H^1\mathcal{V}^{\bm{Y}}_T(\bar{M})$ our considerations so far will imply the claim.
\newline
"$\subseteq$": Let $\bm{X}\in L^2\mathcal{V}^{\bm{Y}}_n(\bar{M})$. By definition we may approximate it by a sequence $\bm{X}_k$ in $\mathcal{V}^{\bm{Y}}_n(\bar{M})$ in the $L^2$-norm. Note that if we can show that $\text{curl}^{-1}(X_k)\in \mathcal{V}^{\bm{Y}}_T(\bar{M})$ for every $k$, then the continuity of this operator will imply $\text{curl}^{-1}(\bm{X})\in H^1\mathcal{V}^{\bm{Y}}_T(\bar{M})$. We know already that $\bm{A}_k:=\text{curl}^{-1}(\bm{X}_k)\in \mathcal{V}_T(\bar{M})$ so we need to show that $\bm{A}_k$ commutes with $\bm{Y}$. By \cref{PT1L1} we have the identity $\text{grad}(g(\bm{A}_k,\bm{Y}))=[\bm{Y},\bm{A}_k]+\bm{Y}\times \bm{X}_k$. Since $\bm{Y}$ and $\bm{X}_k$ are divergence-free we have $\text{curl}(\bm{Y}\times \bm{X}_k)=[\bm{X}_k,\bm{Y}]=0$ since the $\bm{X}_k$ all commute with $\bm{Y}$. Thus we have $\text{curl}([\bm{Y},\bm{A}_k])=0$. On the other hand the $\bm{A}_k$ are also divergence-free and so $[\bm{Y},\bm{A}_k]=\text{curl}(\bm{A}_k\times \bm{Y})$ must be divergence-free. We conclude that $[\bm{Y},\bm{A}_k]\in \mathcal{H}(\bar{M})$ and admits a vector potential. We now recall that $[\bm{Y},\bm{A}_k]=\text{grad}(g(\bm{A}_k,\bm{Y}))+\bm{X}_k\times \bm{Y}$ and note that since $\bm{X}_k$ and $\bm{Y}$ are both tangent, their cross product must be normal to the boundary. Similarly we observe that $\bm{A}_k$ is normal to the boundary, while $\bm{Y}$ is tangent to it, hence $g(\bm{A}_k,\bm{Y})$ vanishes on the boundary and consequently the tangent part of its gradient vanishes. Thus $\text{grad}(g(\bm{A}_k,\bm{Y}))$ is normal to the boundary. We conclude that $[\bm{Y},\bm{A}_k]$ is normal to the boundary, is solenoidal, irrotational and admits a vector potential. It then follows from integration by parts that $[\bm{Y},\bm{A}_k]$ must be the zero vector field as desired.
\newline
"$\supseteq$": Fix any $\bm{A}\in H^1\mathcal{V}^{\bm{Y}}_T(\bar{M})$ and approximate it by a sequence $\bm{A}_k$ in $\mathcal{V}^{\bm{Y}}_T(\bar{M})$ in $H^1$-norm. With a similar reasoning as before it is enough to show that $\bm{X}_k:=\text{curl}(\bm{A}_k)$ commutes with $\bm{Y}$. But this is easy to see. From \cref{PT1L1} we have $\text{grad}(g(\bm{A}_k,\bm{Y}))=\bm{Y}\times \bm{X}_k+[\bm{Y},\bm{A}_k]=\bm{Y}\times \bm{X}_k$ since each $\bm{A}_k$ commutes with $\bm{Y}$. Now we can apply the curl to this equation and keeping in mind that $\bm{Y}$ and the $\bm{X}_k$ are divergence-free we find $[\bm{X}_k,\bm{Y}]=\text{curl}(\bm{Y}\times \bm{X}_k)=0$ as desired. $\square$
\newline
\newline
Observe that we have in particular shown
\begin{equation}
\label{PT16}
\text{curl}^{-1}\left(\mathcal{V}^{\bm{Y}}_n(\bar{M}) \right)=\mathcal{V}^{\bm{Y}}_T(\bar{M}).
\end{equation}
Note further that $L^2\mathcal{V}^{\bm{Y}}_n(\bar{M})$ is an $L^2$-closed subspace of $L^2\mathcal{V}(\bar{M})$ and thus we may define the $L^2$-orthogonal projection operator $\pi:L^2\mathcal{V}(\bar{M})\rightarrow L^2\mathcal{V}^{\bm{Y}}_n(\bar{M})$.
\begin{lem}
\label{PT1L5}
Let $(\bar{M},g)$ be a compact $3$-manifold and $\bm{Y}\in \mathcal{K}(\bar{M})$. Then the operator
\begin{equation}
\label{PT17}
\pi\circ \text{curl}^{-1}: \left(L^2\mathcal{V}^{\bm{Y}}_n(\bar{M}),\norm{\cdot}_{L^2}\right)\rightarrow \left(L^2\mathcal{V}^{\bm{Y}}_n(\bar{M}),\norm{\cdot}_{L^2}\right)
\end{equation}
is a compact, self-adjoint operator. Further, for every $\bm{X}\in L^2\mathcal{V}^{\bm{Y}}_n(\bar{M})$ we have $(\pi\circ \text{curl}^{-1})(\bm{X})\in L^2\mathcal{V}^{\bm{Y}}_n(\bar{M})\cap H^1\mathcal{V}(\bar{M})$ and 
\begin{equation}
\label{PT18}
\text{curl}\left((\pi\circ \text{curl}^{-1})(\bm{X}) \right)=\bm{X}.
\end{equation}
\end{lem}
\underline{Proof of \cref{PT1L5}:} The compactness is immediate from \cref{PT1L4} and the $L^2$-continuity of $\pi$. As for the self-adjointness we may, by an approximation argument, without loss of generality assume that we are given two vector fields $\bm{X},\bm{Z}\in \mathcal{V}^{\bm{Y}}_n(\bar{M})$. Note that $\text{curl}^{-1}(\bm{X})-(\pi\circ \text{curl}^{-1})(\bm{X})$ is $L^2$-orthogonal to $L^2\mathcal{V}^{\bm{Y}}_n(\bar{M})$ by definition of the projection and thus
\[
\langle (\pi\circ \text{curl}^{-1})(\bm{X}),\bm{Z} \rangle_{L^2}=\langle \text{curl}^{-1}(\bm{X}),\bm{Z} \rangle_{L^2}=\langle \text{curl}^{-1}(\bm{X}),\text{curl}\left(\text{curl}^{-1}(\bm{Z})\right)\rangle_{L^2}
\]
\[
=\langle \text{curl}\left(\text{curl}^{-1}(\bm{X})\right), \text{curl}^{-1}(\bm{Z})\rangle_{L^2}= \langle \bm{X},\text{curl}^{-1}(\bm{Z})\rangle_{L^2}=\langle \bm{X},(\pi\circ\text{curl}^{-1})(\bm{Z})\rangle_{L^2},
\]
where we used the integration by parts formula and the boundary conditions of elements in $\mathcal{V}_T(\bar{M})$. This proves self-adjointness.
\newline
For the remaining claims fix any $\bm{X}\in L^2\mathcal{V}^{\bm{Y}}_n(\bar{M})$ and let $\bm{X}_k$ be an $L^2$-approximating sequence in $\mathcal{V}^{\bm{Y}}_n(\bar{M})$. By continuity of $\text{curl}^{-1}$ we know that the sequence $\bm{A}_k:=\text{curl}^{-1}(\bm{X}_k)$ converges in $H^1$ to $\text{curl}^{-1}(\bm{X})$. Further the Hodge-Morrey decomposition allows us to write $\bm{A}_k=\bm{B}_k+\bm{\Gamma}_k$ for suitable $\bm{B}_k\in \mathcal{V}_n(\bar{M})$ and $\bm{\Gamma}_k\in \mathcal{H}(\bar{M})$ since the $\bm{A}_k$ are divergence-free. It is easy to check that the spaces $\mathcal{H}(\bar{M})$ and $\mathcal{V}_n(\bar{M})$ are $L^2$-orthogonal and consequently we have $\pi(\bm{\Gamma}_k)=0$ for all $k$ and hence by continuity
\begin{equation}
\label{PT19}
(\pi\circ \text{curl}^{-1})(\bm{X})=\lim_{k\rightarrow \infty}\pi(\bm{B}_k)\text{ in }L^2.
\end{equation}
We know already that $\bm{B}_k\in \mathcal{V}_n(\bar{M})$ for all $k$. We claim that in fact $[\bm{B}_k,\bm{Y}]\equiv 0$ for all $k$. To see this we recall that according to (\ref{PT16}) we have $[\bm{A}_k,\bm{Y}]\equiv 0$ for all $k$ and since the $\bm{\Gamma}_k$ are irrotational we have $\text{curl}(\bm{B}_k)=\text{curl}(\bm{A}_k)=\bm{X}_k$. Hence we find by \cref{PT1L1}
\begin{equation}
\label{PT110}
\text{grad}(g(\bm{B}_k,\bm{Y}))=[\bm{Y},\bm{B}_k]+\bm{Y}\times \bm{X}_k
\end{equation}
\[
\Rightarrow \text{curl}([\bm{Y},\bm{B}_k])=[\bm{Y},\bm{X}_k]=0,
\]
where we used our standard calculus identity, since $\bm{Y}$ and $\bm{X}_k$ are divergence-free, and that the $\bm{X}_k$ commute with $\bm{Y}$ by choice of our approximating sequence. Once again making use of the fact that $\bm{Y}$ as well as $\bm{B}_k$ are divergence-free we find $[\bm{Y},\bm{B}_k]=\text{curl}(\bm{B}_k\times \bm{Y})$ and thus overall $[\bm{Y},\bm{B}_k]\in \mathcal{H}(\bar{M})$. On the other hand it follows from \cref{PT1L1} and the fact that the $\bm{A}_k$ commute with $\bm{Y}$ that
\[
\text{grad}(g(\bm{A}_k,\bm{Y}))=\bm{Y}\times \bm{X}_k+[\bm{Y},\bm{A}_k]=\bm{Y}\times \bm{X}_k.
\]
Thus $\bm{Y}\times \bm{X}_k$ is a gradient field which in combination with (\ref{PT110}) implies that $[\bm{Y},\bm{B}_k]$ is a gradient field. Now we observe that $\bm{Y}$ as well as $\bm{B}_k$ are both tangent to the boundary and thus so must be their Lie bracket. We overall see that $[\bm{Y},\bm{B}_k]$ is solenoidal, irrotational, tangent to the boundary and admits a scalar potential, thus it must be the zero vector field, i.e. $[\bm{Y},\bm{B}_k]\equiv 0$ for all $k$. From this we obtain $\bm{B}_k\in \mathcal{V}^{\bm{Y}}_n(\bar{M})$ and consequently $\pi(\bm{B}_k)=\bm{B}_k$. Thus (\ref{PT19}) implies that $\bm{B}_k$ converges in $L^2$ to $(\pi\circ \text{curl}^{-1})(\bm{X})$. It now follows from elliptic estimates, \cite[Lemma 2.4.10]{S95}, and the fact that the $\bm{X}_k$ form an $L^2$-Cauchy sequence, that the $\bm{B}_k$ form an $H^1$-Cauchy sequence. Hence they converge in $H^1$ to $(\pi\circ\text{curl}^{-1})(\bm{X})$ which proves the regularity assertion. Lastly, since $\bm{B}_k$ converges to $(\pi\circ\text{curl}^{-1})(\bm{X})$ in $H^1$, their curls converge to $\text{curl}\left((\pi\circ\text{curl}^{-1})(\bm{X})\right)$ in $L^2$. But we recall that we have $\text{curl}(\bm{B}_k)=\text{curl}(\bm{A}_k)=\bm{X}_k$, which converges to $\bm{X}$ in $L^2$. Thus the claim follows. $\square$
\newline
\newline
Observe that we have in particular shown the following: For any fixed $\bm{Y}\in \mathcal{K}(\bar{M})$:
\begin{equation}
\label{PT111}
\forall \bm{X}\in \mathcal{V}^{\bm{Y}}_n(\bar{M}):\text{ }\bm{Y}\times \bm{X}=\text{grad}(g(\text{curl}^{-1}(\bm{X}),\bm{Y})).
\end{equation}
We can now exploit the fact that the constructed operator is compact and self-adjoint
\begin{lem}
\label{PT1L6}
Let $(\bar{M},g)$ be a compact $3$-manifold and $\bm{Y}\in \mathcal{K}(\bar{M})$ with $\mathcal{V}^{\bm{Y}}_n(\bar{M})\neq \{0\}$, then there exists some non-zero $\bm{X}\in L^2\mathcal{V}^{\bm{Y}}_n(\bar{M})\cap H^1\mathcal{V}(\bar{M})$ and a constant $\lambda\in \mathbb{R}\setminus \{0\}$ such that $\text{curl}(\bm{X})=\lambda \bm{X}$.
\end{lem}
\underline{Proof of \cref{PT1L6}:} According to \cref{PT1L5} we know that the operator $\pi\circ \text{curl}^{-1}$ is compact and self-adjoint. Thus there exists a (at most countable) collection of eigenfields $\bm{E}_k$ and corresponding eigenvalues $\mu_k\in \mathbb{R}\setminus \{0\}$ of this operator such that $L^2\mathcal{V}^{\bm{Y}}_n(\bar{M})$ is the $L^2$-orthogonal direct sum of the kernel of this operator and the closure of the span of these eigenfields. Assume for the moment that $L^2\mathcal{V}^{\bm{Y}}_n(\bar{M})=\text{Ker}\left(\pi\circ \text{curl}^{-1}\right)$. Then for every $\bm{X}\in L^2\mathcal{V}^{\bm{Y}}_n(\bar{M})$ we have $(\pi\circ \text{curl}^{-1})(\bm{X})=0$. But then (\ref{PT18}) implies that for each such $\bm{X}$ we have the identity $\bm{X}=\text{curl}\left((\pi\circ \text{curl}^{-1})(\bm{X})\right)=\text{curl}(0)=0$, i.e. $L^2\mathcal{V}^{\bm{Y}}_n(\bar{M})=\{0\}$, which contradicts our assumption $\mathcal{V}^{\bm{Y}}_n(\bar{M})\neq \{0\}$. Thus the operator $(\pi\circ \text{curl}^{-1})$ admits at least one (non-zero) eigenfield $\bm{X}$ corresponding to a non-zero eigenvalue $\mu$. Due to the regularity result of \cref{PT1L5} we have $\bm{X}=\frac{1}{\mu}(\pi\circ\text{curl}^{-1})(\bm{X})\in H^1\mathcal{V}(\bar{M})$ and applying the curl on both sides in combination with (\ref{PT18}) yields $\text{curl}(\bm{X})=\frac{1}{\mu}\bm{X}$. $\square$
\newline
\newline
In order to conclude the proof of \cref{MT1} we are left with establishing the regularity of $\bm{X}$
\begin{lem}
\label{PT1L7}
Let $(\bar{M},g)$ be a compact $3$-manifold and $\bm{Y}\in \mathcal{K}(\bar{M})$. Suppose $\bm{X}\in L^2\mathcal{V}^{\bm{Y}}_n(\bar{M})\cap H^1\mathcal{V}(\bar{M})$ satisfies $\text{curl}(\bm{X})=\lambda \bm{X}$ for some constant $\lambda\neq 0$. Then $\bm{X}\in \mathcal{V}^{\bm{Y}}_n(\bar{M})$.
\end{lem}
\underline{Proof of \cref{PT1L7}:} In order to see that $\bm{X}$ must be smooth, one readily checks that since $\bm{X}\in L^2\mathcal{V}_n(\bar{M})$, its potential $\text{curl}^{-1}(\bm{X})$ (modulo a constant) coincides with the Dirichlet potential of $\bm{X}$, \cite[Theorem 2.2.4]{S95}. The smoothness then follows from the regularity theory of the Dirichlet potential, \cite[Theorem 2.2.6]{S95} and standard Sobolev embeddings \cite[Theorem 1.3.6]{S95}. Thus $\bm{X}\in \mathcal{V}(\bar{M})$. Now since we can approximate $\bm{X}$ in $L^2$ by a sequence $(\bm{X}_k)_k\subset \mathcal{V}_n(\bar{M})$ it follows from the $L^2$-orthogonality of the Hodge-Morrey decomposition, \cite[Theorem 2.4.2, Corollary 3.5.2]{S95}, that in fact $\bm{X}\in \mathcal{V}_n(\bar{M})$. Thus we are left with proving that $[\bm{X},\bm{Y}]\equiv 0$. Observe that $[\bm{X},\bm{Y}]$ contains derivatives acting upon $\bm{X}$ and therefore the fact that we may approximate $\bm{X}$ by a sequence of smooth vector fields commuting with $\bm{Y}$ in $L^2$ does not immediately imply that $\bm{X}$ commutes with $\bm{Y}$. To bypass this problem we note the following: Let $\bm{Z}\in \mathcal{V}_n(\bar{M})$ be any fixed vector field and let for notational simplicity $\bm{W}:=[\bm{X},\bm{Y}]$ which we know is a smooth vector field, then we have by our standard vector calculus identity $[\bm{Z},\bm{Y}]=\text{curl}(\bm{Y}\times \bm{Z})$ and thus
\begin{equation}
\label{PT112}
\langle \bm{W},[\bm{Z},\bm{Y}]\rangle_{L^2}=\langle \bm{W},\text{curl}(\bm{Y}\times \bm{Z})\rangle_{L^2}=\langle \text{curl}(\bm{W}),\bm{Y}\times \bm{Z}\rangle_{L^2}=\langle \text{curl}(\bm{W})\times \bm{Y},\bm{Z}\rangle_{L^2},
\end{equation}
where we used the integration by parts formula, the fact that $\bm{Y}$ and $\bm{Z}$ are both tangent and hence their cross product is normal to the boundary and the scalar triple product rule. However the right hand side of (\ref{PT112}) depends $L^2$-continuously on $\bm{Z}$. Thus if we let $\bm{X}_k$ be a sequence in $\mathcal{V}^{\bm{Y}}_n(\bar{M})$ approximating $\bm{X}$ in $L^2$, we have
\[
\langle [\bm{X},\bm{Y}],[\bm{X},\bm{Y}]\rangle_{L^2}=\langle \bm{W},[\bm{X},\bm{Y}]\rangle_{L^2}=\langle \text{curl}(\bm{W})\times \bm{Y},\bm{X}\rangle_{L^2}
\]
\[
=\lim_{k\rightarrow \infty}\langle \text{curl}(\bm{W})\times \bm{Y},\bm{X}_k\rangle_{L^2}=\lim_{k\rightarrow \infty}\langle \bm{W},[\bm{X}_k,\bm{Y}]\rangle_{L^2}=0
\]
since $[\bm{X}_k,\bm{Y}]\equiv 0$ for each $k$. We conclude $[\bm{X},\bm{Y}]\equiv 0$ as desired. $\square$
\section{Proof of \cref{ExtraP3}, \cref{Extrap4} and \cref{MCA1}}
\underline{Proof of \cref{ExtraP3}:} Assume first $\bm{Y}\equiv 0$. In this case $\mathcal{V}^{\bm{Y}}_n(\bar{M})=\mathcal{V}_n(\bar{M})$. We can then fix any point $p\in M:=\text{int}(\bar{M})$ and an interior chart around that point. Fix any bump function $\rho$ supported in the domain of the chosen chart and set $\bm{X}:=\text{curl}(\rho \bm{B})$, where $\bm{B}$ is any smooth vector field defined in local coordinates on the domain of the chart, such that $\bm{X}$ is not the zero vector field. Then obviously $0\not\equiv \bm{X}\in \mathcal{V}_n(\bar{M})$.
\newline
Now assume $\partial\bar{M}\neq \emptyset$ and suppose we are given any $\bm{Y}\in \mathcal{K}(\bar{M})$. Fix any connected component, $\bar{M}_C$, of $\bar{M}$ which has non-empty boundary. On this component we can write, according to \cref{ExtraT1C4}, $\bm{Y}=\text{curl}(\bm{B})$ for some suitable $\bm{B}\in \mathcal{V}_n(\bar{M}_C)$. If $\bm{Y}$ vanishes everywhere on $\bar{M}_C$, we may define a vector field $\bm{X}$ by setting it to zero outside the component and use the result of the first part to obtain a non trivial vector field, commuting with $\bm{Y}$. Thus suppose now that $\bm{Y}$ does not vanish identically on $\bar{M}_C$. Then we must have $\bm{B}\not\equiv 0$ and we can proceed as usual. By \cref{PT1L1} we have the identity $\text{grad}(g(\bm{Y},\bm{B}))=[\bm{Y},\bm{B}]+\bm{Y}\times \text{curl}(\bm{B})=[\bm{Y},\bm{B}]$, since $\bm{B}$ is a vector potential of $\bm{Y}$. Further, since $\bm{B}\in \mathcal{V}_n(\bar{M}_C)$, it is divergence-free, so that we find $[\bm{Y},\bm{B}]=\text{curl}(\bm{B}\times \bm{Y})$, implying that $[\bm{Y},\bm{B}]\in \mathcal{H}(\bar{M}_C)$. Observe that $\bm{B}$ as well as $\bm{Y}$ are tangent to the boundary and thus so must be $[\bm{Y},\bm{B}]$. Hence $[\bm{Y},\bm{B}]\in \mathcal{H}_N(\bar{M}_C)$ and as we have seen $[\bm{Y},\bm{B}]$ is a gradient field, which implies $[\bm{Y},\bm{B}]\equiv 0$, i.e. $0\not\equiv \bm{B}\in \mathcal{V}^{\bm{Y}}_n(\bar{M}_C)$. Setting $\bm{B}$ to zero outside the connected component gives us the desired vector field.
\newline
Lastly suppose we are given some $\bm{Y}\in \mathcal{K}(\bar{M})$ and that there exists some $\tilde{\bm{Y}}\in \mathcal{K}(\bar{M})$ with $\text{curl}(\tilde{\bm{Y}})\not\equiv 0$. We distinguish two cases. First assume $\text{curl}(\bm{Y})\not\equiv 0$. In that case we can write, according to \cref{PT1C3}, $\bm{Y}=\bm{X}+\bm{\Gamma}$ with suitable $\bm{X}\in \mathcal{V}^{\bm{Y}}_n(\bar{M})$ and $\bm{\Gamma}\in \mathcal{H}_N(\bar{M})$. Since $\bm{Y}$ is not irrotational, but $\bm{\Gamma}$ is, we must have $0\not\equiv \bm{X}\in \mathcal{V}^{\bm{Y}}_n(\bar{M})$. The second case to consider is $\text{curl}(\bm{Y})\equiv 0$. Note that $\bm{Y}$ is tangent to the boundary and divergence-free, i.e. we have $\bm{Y}\in \mathcal{H}_N(\bar{M})$ in this case. Similarly as before we can fix $\tilde{\bm{Y}}$ and decompose it as $\tilde{\bm{Y}}=\bm{X}+\bm{\Gamma}$ with $\bm{X}\in \mathcal{V}_n(\bar{M})\setminus \{0\}$ and $\bm{\Gamma}\in \mathcal{H}_N(\bar{M})$. Since $\bm{Y}\in \mathcal{K}(\bar{M})$ it follows from \cref{PT1L2} that $[\bm{Y},\bm{\Gamma}]\equiv 0$. On the other hand $\tilde{\bm{Y}}\in \mathcal{K}(\bar{M})$ and as we have pointed out already $\bm{Y}\in \mathcal{H}_N(\bar{M})$. Thus the same lemma implies $[\tilde{\bm{Y}},\bm{Y}]\equiv 0$. We obtain $[\bm{Y},\bm{X}]=[\bm{Y},\tilde{\bm{Y}}]-[\bm{Y},\bm{\Gamma}]\equiv 0$, i.e. $0\not\equiv \bm{X}\in \mathcal{V}^{\bm{Y}}_n(\bar{M})$. $\square$
\newline
\newline
\underline{Proof of \cref{MCA1}:} Define the vector spaces $V^N(\bar{M}):=\cap_{j=1}^N\mathcal{V}^{\bm{Y}_j}_n(\bar{M})$ and $V^{T}(\bar{M}):=\cap_{j=1}^N\mathcal{V}^{\bm{Y}_j}_T(\bar{M})$. Then replacing the spaces $\mathcal{V}^{\bm{Y}}_n(\bar{M})$ and $\mathcal{V}^{\bm{Y}}_T(\bar{M})$ by the spaces $V^N(\bar{M})$ and $V^T(\bar{M})$ respectively in the preceding arguments, the results of \cref{PT1L4}-\cref{PT1L7} generalise immediately. Thus we are left with proving that $V^N(\bar{M})\neq \{0\}$ under the given assumptions. If $\text{curl}(\bm{Y}_1)\not\equiv 0$ we decompose $\bm{Y}_1$ as $\bm{Y}_1=\bm{X}+\bm{\Gamma}$ with $\bm{\Gamma}\in \mathcal{H}_N(\bar{M})$ and $\bm{X}\in \mathcal{V}_n^{\bm{Y}_1}(\bar{M})$, where $\bm{X}\not\equiv 0$ as $\text{curl}(\bm{Y}_1)\not\equiv 0$. According to \cref{PT1L2} we have $[\bm{Y}_j,\Gamma]\equiv 0$ for all $1\leq j\leq N$ and consequently by our assumption $[\bm{X},\bm{Y}_j]=[\bm{Y}_1,\bm{Y}_j]-[\bm{\Gamma},\bm{Y}_j]\equiv 0$, which concludes the proof in this case. 
\newline
The remaining case assumes $\partial\bar{M}\neq \emptyset$. We first once more fix some connected component of $\bar{M}$ with non-empty boundary. If all Killing fields vanish on this component, we can argue just like in the case $N=1$. By setting the vector field, which we are about to construct, to zero outside the considered connected component we may without loss of generality assume that $\bar{M}$ is connected and that $\bm{Y}_1\not\equiv 0$. According to \cref{ExtraT1C4} we can write $\bm{Y}_1=\text{curl}(\bm{B})$ for a suitable $\bm{B}\in \mathcal{V}_n(\bar{M})\setminus\{0\}$. We claim that for any fixed $1\leq j \leq N$ the vector field $\bm{Y}_j\times \bm{Y}_1$ is a gradient field. To see this observe first that $\text{curl}(\bm{Y}_j\times \bm{Y}_1)=[\bm{Y}_1,\bm{Y}_j]\equiv 0$, where we used that Killing fields are divergence-free and our assumptions. Hence by means of the Hodge-Morrey-Friedrichs decomposition, \cite[Corollary 3.5.2]{S95}, the vector field $\bm{Y}_j\times \bm{Y}_1$ is a gradient field if and only if it is $L^2$-orthogonal to the space $\mathcal{H}_N(\bar{M})$. Using the scalar triple product rule this follows from integration by parts, keeping in mind that the quantity $\bm{\Gamma}\times \bm{Y}_j$ is normal to the boundary for every $\bm{\Gamma}\in \mathcal{H}_N(\bar{M})$. Indeed if $\bm{\Gamma}\in \mathcal{H}_N(\bar{M})$ is any fixed element we have
\[
\langle \bm{Y}_j\times \bm{Y}_1,\bm{\Gamma}\rangle_{L^2}=\langle \overbrace{\bm{Y}_1}^{=\text{curl}(\bm{B})},\bm{\Gamma}\times \bm{Y}_j\rangle_{L^2}=\langle \bm{B},\text{curl}(\bm{\Gamma}\times \bm{Y}_j)\rangle_{L^2}.
\]
Since $\Gamma$ and $\bm{Y}_j$ are divergence-free we have $\text{curl}(\bm{\Gamma}\times \bm{Y}_j)=[\bm{Y}_j,\bm{\Gamma}]\equiv 0$, in view of \cref{PT1L2}. Thus $\bm{Y}_j\times \bm{Y}_1$ is a gradient field as claimed. Lastly we claim that the vector field $\bm{B}$ commutes with every $\bm{Y}_j$. To see this note that by the Killing property we have $\text{grad}(g(\bm{Y}_j,\bm{B}))=\bm{Y}_j\times \text{curl}(\bm{B})+[\bm{Y}_j,\bm{B}]=\bm{Y}_j\times \bm{Y}_1+[\bm{Y}_j,\bm{B}]$ by choice of $\bm{B}$. From our considerations so far we conclude that $[\bm{Y}_j,\bm{B}]$ is a gradient field. On the other hand, since $\bm{B}\in \mathcal{V}_n(\bar{M})$ it is divergence-free and so our standard argument implies that $[\bm{Y}_j,\bm{B}]$ is divergence-free as well. Thus $[\bm{Y}_j,\bm{B}]\in \mathcal{H}(\bar{M})$ is a gradient field for every $1\leq j\leq N$. But each $\bm{Y}_j$ as well as $\bm{B}$ are tangent to the boundary, so that in fact $[\bm{Y}_j,\bm{B}]\in \mathcal{H}_N(\bar{M})$, which in combination with the gradient property implies that $[\bm{Y}_j,\bm{B}]\equiv 0$ for each $1\leq j\leq N$. Further $\bm{B}$ is not identically zero, which proves the claim. $\square$
\newline
\newline
\underline{Proof of \cref{Extrap4}:} First observe that (ii) trivially implies (iii). The converse implication is provided by \cref{MT1}. Thus it is enough to show the following chain of implications $(i) \Rightarrow (ii)\Rightarrow (iv)\Rightarrow (i)$. The implication $(i)\Rightarrow (ii)$ is easy to see. The condition $\bm{Y}\times \bm{X}\not\equiv 0$ guarantees that $\bm{X}$ is not identically zero. Further since $\lambda>0$, the vector field $\bm{X}$ admits a vector potential. In the case of empty boundary, since boundary conditions become empty conditions, this implies $\bm{X}\in \mathcal{V}_n(M)$, which concludes this implication. To see the implication $(ii)\Rightarrow (iv)$ we define the function $f:=g(\bm{Y},\bm{X})$. Observe that by the Killing field property we have $\text{grad}(f)=\text{grad}(g(\bm{Y},\bm{X}))=[\bm{Y},\bm{X}]+\bm{Y}\times \text{curl}(\bm{X})=\lambda \bm{Y}\times \bm{X}$, by properties of $\bm{X}$. In particular $\text{grad}(f)$ is $g$-orthogonal to $\bm{Y}$, which implies $df(\bm{Y})=\bm{Y}(f)=g(\bm{Y},\text{grad}(f))\equiv 0$. Suppose that $df\equiv 0$, then since $\lambda \neq 0$, we must have $\bm{Y}\times \bm{X}\equiv 0$. Note that by assumption $\bm{Y}$ is not the zero vector field and that $\text{curl}(\bm{Y})\equiv 0$ implies, by \cref{MP7}, that $g(\bm{Y},\bm{Y})$ is locally constant. Thus by connectedness $\bm{Y}$ is nowhere vanishing. Since $\bm{X}$ and $\bm{Y}$ are everywhere collinear, it follows that there is some smooth function $h$ on $M$ with $\bm{X}=h\bm{Y}$. We then have the identity
\[
\lambda h\bm{Y}=\lambda \bm{X}=\text{curl}(\bm{X})=\text{curl}(h\bm{Y})=\text{grad}(h)\times \bm{Y}+\text{curl}(\bm{Y})h=\text{grad}(h)\times \bm{Y},
\]
since $\bm{Y}$ is irrotational. Multiplying the equation by $\bm{Y}$ and using the fact that $\bm{Y}$ is no-where vanishing and that $\lambda \neq 0$, we conclude $h\equiv 0$ and consequently $\bm{X}=h\bm{Y}\equiv 0$, which is a contradiction. Thus $df\not\equiv 0$, which concludes this step. The implication $(iv)\Rightarrow (i)$ can be seen as follows. Fix any $f\in C^{\infty}(M)$ which is not constant and satisfies $df(\bm{Y})\equiv 0$. Defining $c:=\frac{1}{\text{vol}(M)}\int_Mf\omega_g$ and $\phi:=f-c$, we see that $\phi$ is a smooth function which integrates to zero, or equivalently it is $L^2$-orthogonal to the space of constant functions. It then follows that $d\phi=df\not\equiv 0$, hence $\phi$ is not constant, and that $d\phi(\bm{Y})=df(\bm{Y})\equiv 0$. This implies that the set $H^1_{T,\bm{Y}}(M)$, as defined in (\ref{AppendixA1}) in \cref{AppendixA} contains a nonzero element. The claim then follows from \cref{AppendixAC3}. $\square$
\section{Proof of \cref{MT3}}
We will first prove a more general result and then show how exactly it applies to symmetric Beltrami fields, \cref{MT3}, as well as to symmetric, real analytic vector fields with non-trivial helicity, see \cref{MR4}.
\begin{prop}
\label{PT3P1}
Let $(\bar{M},g)$ be a compact, connected, real analytic $3$-manifold. Let $\bm{Y},\bm{X}\in \mathcal{V}^{\omega}(\bar{M})$ be real analytic vector fields which are tangent to the boundary, not everywhere collinear and which commute with each other, $[\bm{X},\bm{Y}]\equiv 0$. If there exists a real analytic function $f:\bar{M}\rightarrow \mathbb{R}$ with $\bm{X}\times \bm{Y}=\text{grad}(f)$, then the conclusions (i)-(v) of \cref{MT3} apply to $\bm{X}$ as well as $\bm{Y}$ (in both cases with the same set $\Gamma$ and the same tori decomposition of its complement).
\end{prop}
\underline{Proof of \cref{PT3P1}:} The proof follows the exposition of the proof of Arnold's structure theorem given in \cite[Chapter II Theorem 1.2]{AK98}. Let $K:=\{p\in \bar{M}| df(p)=0 \}$ and $C:=\{p\in \partial\bar{M}| (d\iota^{\#}f)(p)=0 \}$, where $d$ denotes the exterior derivative and $\iota^{\#}$ is the pullback via the inclusion map. Define $\Gamma:=f^{-1}\left(f(K)\cup f(C) \right)$ and observe that since $\bm{X}$ and $\bm{Y}$ are tangent to the boundary by assumption, the gradient of $f$, being equal to their cross product, is normal to the boundary and thus $C=\partial\bar{M}$. Let $M:=\text{int}(\bar{M})$, $f|_M:M\rightarrow \mathbb{R}$ and $K_M:=\{p\in M| df|_M(p)=0\}$, then we have the identity $\Gamma=f^{-1}(f(\partial \bar{M}))\cup f|^{-1}_{M}(f|_M(K_M))$. Since the boundary is compact it has finitely many connected components and since the restriction of $f$ to the boundary is locally constant, there exist constants $c_1,c_2,\dots, c_N\in \mathbb{R}$ for some $N\in \mathbb{N}$ such that $f^{-1}(f(\partial\bar{M}))=\cup_{i=1}^Nf^{-1}(c_i)=\partial\bar{M}\cup \cup_{i=1}^Nf|^{-1}_M(c_i)$. Let $\Gamma_M:=f|^{-1}_{M}(f|_M(K_M))$, then we obtain
\[
\Gamma=\partial\bar{M} \sqcup  \left(\Gamma_M\cup \bigcup_{i=1}^Nf|^{-1}_M(c_i) \right),
\]
where $\sqcup$ indicates that the union is disjoint. We claim that the latter can be written as a countable union of semianalytic sets. Since $f$ is real analytic, it is clear that each of the sets $f|^{-1}_M(c_i)$ are real analytic and hence in particular semianalytic. As for $\Gamma_M$ we note that $M$ itself is not necessarily compact, however for any $p\in K_M$ we can fix a (analytic) chart $\mu$ of $M$ around $p$ and choose $0<r_p$ so small that the closure of $B_{r_p}(\mu(p))$ is still contained in the charts image. Define $B_{r_p}(p):=\mu^{-1}\left(B_{r_p}(\mu(p)) \right)$. This gives rise to an open cover of $K_M$ and by second countability we can extract a countable subcover $B_{r_i}(p_i)$, $i\in \mathbb{N}$. Note that the topological closures of the $B_{r_i}(p_i)$ are semianalytic subsets of $M$ and are compact. Since $K_M$ is also a semianalytic subset of $M$ it follows that $K_i:=\text{clos}(B_{r_i}(p_i))\cap K_M$ is a countable family whose union is $K_M$ and which consists of compact, semianalytic sets. It follows that $f|_M(K_i)$ are subanalytic subsets of $\mathbb{R}$. Since $\mathbb{R}$ is one-dimensional $f|_M(K_i)$ is semianalytic \cite[Theorem 6.1]{BM88} and hence so are the $\Gamma_i:=f|^{-1}_M(f|_M(K_i))$. Since the union of the $K_i$ equals $K_M$, the union of the $\Gamma_i$ equals $\Gamma_M$ which proves that $\Gamma_M$ is the countable union of semianalytic sets. Observe that none of the $f|^{-1}_M(c_i)$ can contain an interior point, since otherwise $f|_M$ would be constant by real analyticity. Similarly by Sard's theorem the set $f|_M(K_i)$ is a null set and thus has no interior points. Note that the $K_i$ are compact and hence $f|_M(K_i)$ are no-where dense subsets. It then follows that the $\Gamma_i$ are no-where dense subsets of $M$, see also \cite[Chapter II Lemma 1.12]{AK98}. We then obtain from \cite[Corollary 2.11]{BM88} that each of the $\Gamma_i$, as well as the $f|^{-1}_M(c_i)$ can be written as disjoint, countable unions of connected real analytic submanifolds of $M$. Since none of these sets has interior points we conclude that $\Gamma\setminus \partial\bar{M}$ is a countable union of disjoint, real analytic submanifolds of $M$ of dimension at most $2$. This shows that $\Gamma$ is an $\mathcal{H}^2$-countably-$2$ rectifiable set with respect to the natural metric $d_g$ induced by the Riemannian metric $g$ on $\bar{M}$. The compactness of $\Gamma$ is immediate from its definition, so that this concludes the characterisation of $\Gamma$.
\newline
As for items (ii) and (iii) of \cref{MT3} we observe that for every $q\in \bar{M}\setminus \Gamma$ we have $f^{-1}(f(q))\subseteq \bar{M}\setminus \Gamma$. Note further that due to this we have for each $q\in \bar{M}\setminus \Gamma\subseteq M$ the equality $f|^{-1}_M(f(q))=f^{-1}(f(q))$ and thus the sets $f|^{-1}_M(f(q))$ are all compact. In addition $f(q)$ are regular values of $f|_M$ for each $q\in \bar{M}\setminus \Gamma$. Thus the sets $ f|^{-1}_M(f(q))\subseteq \bar{M}\setminus \Gamma$ are all compact, real analytically embedded $2$-dimensional submanifolds without boundary of $M$. Further observe that $\text{grad}(f)$ is always normal to its regular level sets and thus, since $\text{grad}(f)=\bm{X}\times \bm{Y}$, both vector fields $\bm{X}$ and $\bm{Y}$ are tangent to the level sets, which implies that the connected components of each such level set are invariant under the flow of both these vector fields. Note that the gradient of $f$ does not vanish on these surfaces and hence $\bm{X}$ and $\bm{Y}$ must be linearly independent on each connected component of these level sets. So if we fix any such connected component we may restrict the vector fields $\bm{X}$ and $\bm{Y}$ to this submanifold and see that they are linearly independent at each point and that their flows commute by assumption. Thus (ii) and (iii) now follows from the construction of angular coordinates as in the case of Liouville's theorem \cite[Chapter 10]{A89}. Item (v) follows immediately since we have already argued that $\bar{M}\setminus \Gamma$ is the union of invariant tori and hence itself invariant under the flow. Lastly for item (iv) we may fix any such invariant torus $T_i$ and observe that $\text{grad}(f)$ never vanishes on $T_i$. Thus we can find a whole open neighbourhood $U$ of $T_i$ contained in $\bar{M}\setminus \Gamma$ on which $\text{grad}(f)$ does not vanish. Then we may consider the vector field $\frac{\text{grad}(f)}{g(\text{grad}(f),\text{grad}(f))}$ and let $\gamma_p$ denote the integral curve of this vector field starting at a given $p\in U$. One can then use similar arguments as in the proof of the product neighbourhood theorem, \cite{M65}, to show that there exists some $\epsilon>0$ and suitable open neighbourhood $W\subseteq \bar{M}\setminus \Gamma$ of $T_i$ such that
\begin{equation}
\label{PT31}
\psi:(-\epsilon,\epsilon)\times T_i\rightarrow W, (t,p)\mapsto \gamma_p(t)
\end{equation}
defines a diffeomorphism. Now note that for fixed $p\in T_i$ we have $\frac{d}{dt}(f\circ \gamma_p)(t)=1$ by definition of the vector field and $\gamma_p$. Thus we have $(f\circ \gamma_p)(t)=t+c_p$ for some constant $c_p\in \mathbb{R}$ which possibly depends on $p$. However we have $c_p=f(\gamma_p(0))=f(p)\in f(T_i)$ and we recall that each $T_i$ is the connected component of a regular level set of $f$ and thus $f(p)$ has the same value for every $p\in T_i$, i.e. the constant $c_p\equiv c$ is independent of $p$. This implies that for fixed $|t|<\epsilon$ we have for each $p\in T_i$ $f(\psi_t(p))=t+c$, i.e. $\psi_t(T_i)\subseteq f^{-1}(t+c)$. Note that $\psi_t$ is an injective, continuous map between $T_i$ and the regular level set $f^{-1}(t+c)$ and hence is an open map. Further $T_i$ is compact and thus the image $\psi_t(T_i)$ is an open, closed and non-empty subset of $f^{-1}(t+c)$ and therefore must coincide with one of its connected components. Thus is an invariant torus itself as claimed. Lastly it is not hard to verify that if we let $\epsilon_1,\epsilon_2>0$ be maximal with the property that $\psi:(-\epsilon_1,\epsilon_2)\times T_i\rightarrow \bar{M}\setminus \Gamma$ is a well-defined diffeomorphism onto an open set, then it gives rise to a diffeomorphism onto the connected component containing $T_i$. $\square$
\newline
\newline
\underline{Proof of \cref{MT3}:} All we need to do is to verify that the requirements of \cref{PT3P1} are satisfied. Recall that in our scenario $\bm{Y}\in \mathcal{K}^{\omega}(\bar{M})$ is a Killing field and $\bm{X}$ is a real analytic Beltrami field corresponding to a non-zero eigenvalue, tangent to the boundary and that these vector fields commute and are not everywhere collinear by assumption. So we only need to show the existence of a real analytic function $f:\bar{M}\rightarrow \mathbb{R}$ such that $\bm{Y}\times \bm{X}=\text{grad}(f)$, then \cref{PT3P1} will apply. But this follows easily from \cref{PT1L1}, since we have
\[
\text{grad}(g(\bm{X},\bm{Y}))=\bm{Y}\times \text{curl}(\bm{X})+[\bm{Y},\bm{X}]=\lambda \bm{Y}\times \bm{X},
\]
since $\bm{X}$ and $\bm{Y}$ commute and $\bm{X}$ is a Beltrami field. Obviously the function $f:=\frac{g(\bm{X},\bm{Y})}{\lambda}$ is real analytic and satisfies the requirements. $\square$
\begin{rem}
\label{PT3R2}
Note that the tori we constructed in \cref{PT3P1} were regular level sets of the function $f$ and thus the proof of \cref{MT3} in particular implies that $\bm{X}$ and $\bm{Y}$ are tangent to the regular level sets of $g(\bm{X},\bm{Y})$. In view of \cref{ME8} this statement generalises \cite[Theorem 9]{C99}, which dealt with the rotationally symmetric situation in Euclidean space.
\end{rem}
\section{Proof of \cref{MTE2}}
\underline{Proof of \cref{MTE2}:} If we are in case (ii), i.e. if $\text{curl}(\bm{Y})$ is a constant multiple of $\bm{Y}$, then the statement about the geodesics follows immediately from \cref{MP7}, whose proof will be given in the next section. The remaining part of (ii) is a consequence of the results presented in \cref{AppendixA}. Namely if $\bm{Y}(f)\equiv 0$ for some non-constant $C^1$- function $f$, then we may set $c:=\frac{1}{\text{vol}(M)}\int_{M}f\omega_g$ and define $f_c:=f-c$. Then $f_c$ is of class $C^1$, integrates to zero, i.e. is $L^2$-orthogonal to the space of constant functions, satisfies $\bm{Y}(f_c)=\bm{Y}(f)\equiv 0$ and is not constant. It then follows from \cref{AppendixAC3} that there exists some eigenfield $\bm{X}$ of curl, corresponding to a non-zero eigenvalue, which commutes with $\bm{Y}$ and is not everywhere collinear with $\bm{Y}$. The assertion then follows from \cref{MT3} and \cref{MR4}. We are left with considering the first case, that is we assume from now on that there does not exist any $\kappa\in \mathbb{R}$ with $\text{curl}(\bm{Y})=\kappa \bm{Y}$. This assumption in particular implies that $\text{curl}(\bm{Y})\not\equiv 0$ and hence \cref{MT1} in combination with \cref{ExtraP3} guarantees us the existence of a real analytic, strong Beltrami field $\bm{X}$ which commutes with $\bm{Y}$. Now if $\bm{Y}$ and $\bm{X}$ are linearly independent at at least one point, then \cref{MT3} in combination with the 4th item of \cref{MR4} imply that the conclusions of the structure theorem apply to $\bm{Y}$. Thus we may now without loss of generality assume that $\bm{Y}$ and $\bm{X}$ are everywhere collinear. It then follows from \cref{PT1L1} that
\[
\text{grad}(g(\bm{Y},\bm{X}))=[\bm{Y},\bm{X}]+\lambda\bm{Y}\times \bm{X}\equiv 0,
\]
where $\lambda$ is the corresponding eigenvalue of $\bm{X}$ and we used that the vector fields commute and are everywhere collinear by assumption. Thus there exists some $c\in \mathbb{R}$ with $g(\bm{Y},\bm{X})=c$. By real analyticity and since $\bm{X}$ is not the zero vector field we know that the set $U:=\{p\in M| \bm{X}(p)\neq 0 \}$ is an open and dense subset of $M$. Since $\bm{X}$ and $\bm{Y}$ are everywhere collinear, there must exist a real analytic function $f:U\rightarrow \mathbb{R}$ with $\bm{Y}=f\bm{X}$ on $U$. In fact we have
\begin{equation}
\label{PTE1}
c=g(\bm{X},\bm{Y})=fg(\bm{X},\bm{X})\Leftrightarrow f=\frac{c}{g(\bm{X},\bm{X})}\text{ on }U,
\end{equation}
where we used that $\bm{X}$ does not vanish on $U$. Now if $c=0$, then $f\equiv 0$ and consequently $\bm{Y}\equiv 0$ on $U$, which by density of $U$ in $M$ implies that $\bm{Y}$ is the zero vector field. This contradicts our assumption $\text{curl}(\bm{Y})\neq \kappa \bm{Y}$. Thus we must have $c\neq 0$. But then we obtain from (\ref{PTE1}) and the relation $\bm{Y}=f\bm{X}$
\[
\bm{X}=\frac{g(\bm{X},\bm{X})}{c}\bm{Y}.
\]
For notational simplicity let us set $F:=\frac{g(\bm{X},\bm{X})}{c}$. Then on the one hand, since $\bm{X}$ and $\bm{Y}$ are divergence-free, we obtain
\begin{equation}
\label{PTE2}
0=\text{div}(\bm{X})=\text{div}(F\bm{Y})=g(\text{grad}(F),\bm{Y})+F\text{div}(\bm{Y})=g(\text{grad}(F),\bm{Y}).
\end{equation}
On the other hand, since $\bm{X}$ is a Beltrami field, we find
\[
\lambda \bm{X}=\text{curl}(\bm{X})= \text{curl}(F\bm{Y})=\text{grad}(F)\times \bm{Y}+F \text{curl}(\bm{Y}).
\]
Now recall that we are in the case where $\bm{X}$ and $\bm{Y}$ are everywhere collinear and so we may take the cross product of this equation by $\bm{Y}$ from the left to obtain
\[
0=F\bm{Y}\times \text{curl}(\bm{Y})+\bm{Y}\times \left(\text{grad}(F)\times \bm{Y} \right).
\]
By the vector triple product rule we have $\bm{Y}\times \left(\text{grad}(F)\times \bm{Y} \right)=\text{grad}(F)g(\bm{Y},\bm{Y})-\bm{Y}g(\text{grad}(F),\bm{Y})$. Note that the last term is zero by (\ref{PTE2}) and thus we overall obtain the identity
\[
0=F\bm{Y}\times \text{curl}(\bm{Y})+g(\bm{Y},\bm{Y})\text{grad}(F)\text{ on }U.
\]
Assume for the time being that $\bm{Y}\times \text{curl}(\bm{Y})\equiv 0$ on $U$, then the above equation gives us $0=g(\bm{Y},\bm{Y})\text{grad}(F)$ on $U$. Recall that $\bm{Y}=\frac{c}{g(\bm{X},\bm{X})}\bm{X}$ for some $c\neq 0$ and that $\bm{X}$ does not vanish on $U$ by definition of $U$. Hence $\bm{Y}$ also does not vanish on $U$ and thus we conclude $0\equiv \text{grad}(F)$ on $U$. But by definition of $F$ this means that $\text{grad}(g(\bm{X},\bm{X}))\equiv 0$ on $U$ and by density of $U$ in $M$ this extends to all of $M$. In other words $g(\bm{X},\bm{X})$ is constant on $M$ and in conclusion $f=\frac{c}{g(\bm{X},\bm{X})}$ is just a non-zero constant. We therefore have $\bm{Y}=f\bm{X}$ on $U$ for some constant $f\neq 0$. Once again a density argument implies that $\bm{Y}$ is just a non-zero, constant multiple of $\bm{X}$ on $M$ and therefore a Beltrami field, which contradicts our assumption $\text{curl}(\bm{Y})\not\equiv \kappa \bm{Y}$ for any constant $\kappa$. We conclude that we have overall shown that if $\bm{Y}$ and $\bm{X}$ are everywhere collinear and $\text{curl}(\bm{Y})$ is not a constant multiple of $\bm{Y}$, then $\bm{Y}$ and $\text{curl}(\bm{Y})$ are at at least one point linearly independent. Now we observe that by \cref{PT1L1} we have the identity
\begin{equation}
\label{PTE3}
\text{grad}(g(\bm{Y},\bm{Y}))=[\bm{Y},\bm{Y}]+\bm{Y}\times \text{curl}(\bm{Y})=\bm{Y}\times \text{curl}(\bm{Y}).
\end{equation}
Since $\bm{Y}$ and $\text{curl}(\bm{Y})$ are both divergence-free we may take the curl of this equation and obtain from our standard vector calculus identity $[\bm{Y},\text{curl}(\bm{Y})]\equiv 0$. Thus $\bm{Y}$ and $\text{curl}(\bm{Y})$ satisfy the requirements of \cref{PT3P1} (here the tangent to the boundary condition is an empty condition since the boundary is empty). Thus the conclusions of the structure theorem apply to $\bm{Y}$ as claimed.
\newline
Observe lastly that the set $K_M$, which we defined in the proof of \cref{PT3P1}, is compact if the boundary is empty, since it is a closed subset of $M$. But then we can not only find a countable cover of $K_M$ by the $K_i$, but in fact a finite cover. In conclusion $\Gamma_M=\Gamma$ (in the case of empty boundary) is a finite union of $\Gamma_i$'s, each of which is semianalytic. Thus $\Gamma$ is semianalytic as a finite union of such sets and $M\setminus \Gamma$ is semianalytic as a complement of such a set. Finally both of these sets are relatively compact subsets of $M$, because $M$ is compact. It follows from \cite[Corollary 2.7]{BM88} that both sets have at most finitely many connected components. $\square$
\section{Proof of remaining claims}
\underline{Proof of \cref{ML5}:} It follows from $\bm{Y}\times \bm{X}\equiv 0$ in particular that $\text{curl}(\bm{Y}\times \bm{X})\equiv 0$ and since both vector fields are divergence-free this in particular implies that they commute. Thus $\bm{X}\in \mathcal{V}^{\bm{Y}}_n(\bar{M})$ and by (\ref{PT16}) we find $\bm{A}:=\text{curl}^{-1}(\bm{X})\in \mathcal{V}^{\bm{Y}}_T(\bar{M})$. Now it follows from \cref{PT1L1} that 
\[
\text{grad}(g(\bm{A},\bm{Y}))=\bm{Y}\times \bm{X}+[\bm{Y},\bm{A}]=\bm{Y}\times \bm{X}\equiv 0,
\]
the latter by assumption. By connectedness $g(\bm{A},\bm{Y})$ is constant and since the boundary is non-empty this constant must be zero, because $\bm{A}$ is normal to the boundary while $\bm{Y}$ is tangent to it. We conclude $g(\bm{A},\bm{Y})\equiv 0$. Now since $\bm{X}$ is a Beltrami field it follows from a unique continuation result, \cite{AKS62}, that the set $U$ of points at which $\bm{X}$ does not vanish is an open and dense subset of $\bar{M}$. Since $\bm{Y}\times \bm{X}\equiv 0$ these vector fields must be everywhere linearly dependent and since $\bm{X}$ does not vanish on $U$ there exists a smooth function $f:U\rightarrow \mathbb{R}$ with $\bm{Y}=f\bm{X}$. Now since $\bm{X}\in \mathcal{V}^{\bm{Y}}_n(\bar{M})$ and $\bm{X}$ is a Beltrami field we have by \cref{PT1L1}
\[
\text{grad}(g(\bm{X},\bm{Y}))=\lambda \bm{Y}\times \bm{X}+[\bm{Y},\bm{X}]\equiv 0
\]
just like before. Thus there exists some constant $c\in \mathbb{R}$ with $g(\bm{Y},\bm{X})\equiv c$. In particular we have the following identity on $U$
\[
c=fg(\bm{X},\bm{X})\Leftrightarrow f=\frac{c}{g(\bm{X},\bm{X})}.
\]
If $c=0$, this implies $f\equiv 0$ and consequently that $\bm{Y}$ is identical zero on $U$ which is a dense subset of $\bar{M}$, i.e. $\bm{Y}$ is the zero vector field. So assume now that $c\neq 0$, then $f$ is no-where vanishing on $U$ and consequently we find $0=g(\bm{A},\bm{Y})=fg(\bm{A},\bm{X})\Rightarrow 0=g(\bm{A},\bm{X})$ on $U$. Again by a density argument this holds on all of $\bar{M}$. Then since $\bm{X}$ is a Beltrami field this implies
\[
0=\lambda\langle \bm{A},\bm{X}\rangle_{L^2}=\langle \bm{A},\text{curl}(\bm{X})\rangle_{L^2}=\langle \text{curl}(\bm{A}),\bm{X}\rangle_{L^2}=\langle \bm{X},\bm{X}\rangle_{L^2},
\]
since $\bm{A}$ is a vector potential of $\bm{X}$ and normal to the boundary. Thus we must have in this case $\bm{X}\equiv 0$ which contradicts the fact that $\bm{X}$ is an eigenfield. $\square$
\newline
\newline
\underline{Proof of \cref{MP7}:} According to \cref{PT1L1} we have
\[
\text{grad}(g(\bm{Y},\bm{Y}))=\bm{Y}\times \text{curl}(\bm{Y})+[\bm{Y},\bm{Y}]=\bm{Y}\times \text{curl}(\bm{Y}).
\]
Now if $\bm{Y}$ has a vanishing curl or is a Beltrami field the term $\bm{Y}\times \text{curl}(\bm{Y})$ vanishes and the first item follows. Using the following vector calculus identity $\nabla_{\bm{Y}}\bm{Y}=\frac{1}{2}\text{grad}(g(\bm{Y},\bm{Y}))-\bm{Y}\times \text{curl}(\bm{Y})$, we see that $\nabla_{\bm{Y}}\bm{Y}$ must vanish as well.
\newline
As for the second part we observe that if the curl of $\bm{Y}$ vanishes we have by the connectedness assumption and from the first item that $g(\bm{Y},\bm{Y})=c$ for some constant $c\in \mathbb{R}$. If we now restrict $\bm{Y}$ to the compact, connected component of the boundary which is not homeomorphic to the $2$-torus, we see that $\bm{Y}$ (being tangent to the boundary) must vanish at at least one point by virtue of the Poincar\'{e}-Hopf theorem and thus $c=0$. We conclude that $\bm{Y}$ is identical zero in this case.
\newline
Now assume that $\bm{X}\in \mathcal{V}(\bar{M})$ is a Beltrami field such that $\bm{Y}$ and $\bm{X}$ are everywhere collinear. It follows again from \cref{PT1L1} that
\[
\text{grad}(g(\bm{Y},\bm{X}))=[\bm{Y},\bm{X}]+\lambda \bm{Y}\times \bm{X}=0,
\]
where $\lambda$ is the corresponding eigenvalue, where we used that $\bm{Y}$ and $\bm{X}$ are everywhere collinear and that $\text{curl}(\bm{Y}\times \bm{X})=[\bm{Y},\bm{X}]$ since both vector fields are divergence-free. Similarly there exists a constant $c\in \mathbb{R}$ with $g(\bm{Y},\bm{X})=c$ and once again restricting $\bm{Y}$ to the connected component we see that $\bm{Y}$ must vanish at at least one point and thus $c=0$. We conclude $g(\bm{X},\bm{Y})\equiv 0$. Now we can argue just like in the proof of \cref{ML5} that the set of points $U$ at which $\bm{X}$ does not vanish is an open and dense subset of $\bar{M}$ and hence we can find a smooth function $f:U\rightarrow \mathbb{R}$ with $\bm{Y}=f\bm{X}$. But then we have $0=g(\bm{Y},\bm{X})=fg(\bm{X},\bm{X})$ and thus $f\equiv 0$ since $\bm{X}$ is no-where vanishing on $U$. Hence $\bm{Y}$ must vanish on a dense subset which means it is the zero vector field. $\square$
\newline
\newline
So far we did not yet prove the last item of \cref{MR4}. Let us formulate the precise statement as a proposition
\begin{prop}
\label{PRP1}
Let $(\bar{M},g)$ be a compact, connected, real analytic $3$-manifold with non-empty boundary, $\bm{Y}\in \mathcal{K}^{\omega}(\bar{M})\setminus \{0\}$ and $\bm{X}\in \mathcal{V}^{\bm{Y}}_n(\bar{M})$ be a real analytic vector field which satisfies $ \mathcal{H}(\bm{X}):=\langle \bm{X},\text{curl}^{-1}(\bm{X})\rangle_{L^2}\neq 0$, then the conclusions of the structure theorem apply to $\bm{X}$.
\end{prop}
\underline{Proof of \cref{PRP1}:} We want to apply \cref{PT3P1}. Note that $\bm{Y}$ and $\bm{X}$ are tangent to the boundary and that by the definition of the space $\mathcal{V}^{\bm{Y}}_n(\bar{M})$ both vector fields commute. We first claim that they are not everywhere collinear. Assume the opposite, namely that $\bm{Y}\times \bm{X}\equiv 0$. Then just like in the proof of \cref{ML5} we obtain $g(\bm{A},\bm{Y})\equiv 0$, where $\bm{A}:=\text{curl}^{-1}(\bm{X})$. Since $\bm{Y}$ is real analytic and it is not the zero vector field, the set of points $U$ at which $\bm{Y}$ does not vanish must be an open and dense subset of $\bar{M}$. Hence we can find a function $F:U\rightarrow \mathbb{R}$ with $\bm{X}=F\bm{Y}$ and conclude $g(\bm{A},\bm{X})=Fg(\bm{A},\bm{Y})\equiv 0$. Therefore $0=\langle \bm{A},\bm{X}\rangle_{L^2}=\mathcal{H}(\bm{X})$ which contradicts our assumption. Thus $\bm{X}$ and $\bm{Y}$ are at at least one point linearly independent. We are left with showing that there exists a real analytic function $f:\bar{M}\rightarrow \mathbb{R}$ with $\bm{Y}\times \bm{X}=\text{grad}(f)$. But recall that according to (\ref{PT111}) we have
\[
\bm{Y}\times \bm{X}=\text{grad}(g(\bm{A},\bm{Y})).
\]
So we can choose $f:=g(\bm{A},\bm{Y})$. Note that the standard Hodge-Morrey decomposition a priori only tells us that $\bm{A}$ is smooth and so we know that $f$ is smooth. To see that it must in fact be real analytic observe that $\text{grad}(f)=\bm{Y}\times \bm{X}$ is real analytic since both vector fields are by assumption. But it is not hard to see that a smooth function with a real analytic gradient is real analytic. $\square$
\appendix
\section{Symmetry constrained elliptic scalar equations and symmetric Beltrami flows}
\label{AppendixA}
In this section we study a symmetry constrained Laplacian Dirichlet eigenfunction problem and will explain how this scalar problem relates to the existence of symmetric Beltrami fields, provided, the involved Killing field is of Beltrami type. To this end we recall that if $(\bar{M},g)$ is a compact, connected $3$-manifold, we can define the $H^1$-norm on $C^{\infty}(\bar{M})$ by $\norm{f}^2_{H^1}:=\norm{f}^2_{L^2}+\norm{df}^2_{L^2}$, where the latter are the corresponding $L^2$-norms on functions and $1$-forms induced by the metric $g$. We define $H^1(\bar{M})$ to be the completion of $C^{\infty}(\bar{M})$ with respect to the norm $\norm{\cdot}_{H^1}$. In fact the norm is induced by an inner product, so that $H^1(\bar{M})$ is a Hilbert space. Given some $\bm{Y}\in \mathcal{K}(\bar{M})$ we denote by $\omega^1_{\bm{Y}}$ the associated $1$-form via $g$ and by $\langle\cdot,\cdot\rangle_g$ the pointwise fibre product on the cotangent bundle. Further we let $\mathcal{H}^0_D(\bar{M}):=\{f\in C^{\infty}(\bar{M})|f\text{ is constant and }f|_{\partial\bar{M}}\equiv 0 \}$ and denote by $\left(\mathcal{H}^0_D(\bar{M})\right)^{\perp}$ its $L^2$-orthogonal complement. We then define the space
\begin{AppA}
\label{AppendixA1}
H^1_{T,\bm{Y}}(\bar{M}):=\{f\in H^1(\bar{M})|f|_{\partial\bar{M}}\equiv 0,\text{ }\langle df,\omega^1_{\bm{Y}}\rangle_g\equiv 0\text{ and }f\in \left(\mathcal{H}^0_D(\bar{M})\right)^{\perp} \}.
\end{AppA}
Observe that this space is a closed subspace of $H^1(\bar{M})$ for every fixed $\bm{Y}\in \mathcal{K}(\bar{M})$ with respect to the metric $\norm{\cdot}_{H^1}$. Indeed the trace theorem, \cite[Theorem 1.3.7]{S95}, implies that the boundary conditions are preserved under $H^1$-convergence. Also $H^1$-convergence in particular implies $L^2$-convergence, so that the property of $L^2$-orthogonality to the subspace $\mathcal{H}^0_D(\bar{M})$ is also preserved. As for the remaining condition $\langle df,\omega^1_{\bm{Y}}\rangle_g\equiv 0$ we consider first the following map
\[
F_{\bm{Y}}:\left(C^{\infty}(\bar{M}),\norm{\cdot}_{H^1} \right)\rightarrow \left(C^{\infty}(\bar{M}),\norm{\cdot}_{L^2} \right), f\mapsto \langle df,\omega^1_{\bm{Y}}\rangle_g.
\]
Fix $f,h\in C^{\infty}(\bar{M})$ and observe that by compactness of $\bar{M}$ and the Cauchy-Schwarz inequality (applied fibre-wise) we have
\[
\left|\langle df,\omega^1_{\bm{Y}}\rangle_g-\langle dh,\omega^1_{\bm{Y}}\rangle_g\right|^2=\left|\langle df-dh,\omega^1_{\bm{Y}}\rangle_g\right|^2\leq \langle df-dh,df-dh\rangle_g \langle \omega^1_{\bm{Y}},\omega^1_{\bm{Y}}\rangle_g\leq C \langle df-dh,df-dh\rangle_g
\]
for some constant $C>0$ independent of $f$ and $h$. Integrating this inequality yields
\[
\norm{F_{\bm{Y}}(f)-F_{\bm{Y}}(h)}^2_{L^2}\leq C\norm{df-dh}^2_{L^2}\leq C\norm{f-h}^2_{H^1}.
\]
Thus the function $F_{\bm{Y}}$ is a linear bounded map and hence extends uniquely to a continuous linear map defined on $\left(H^1(\bar{M}),\norm{\cdot}_{H^1}\right)\rightarrow \left(L^2(\bar{M}),\norm{\cdot}_{L^2}\right)$, denoted in the same way, where $L^2(\bar{M})$ denotes the $L^2$-completion of $C^{\infty}(\bar{M})$. It is then clear from continuity of this map that the property $0\equiv \langle df,\omega^1_{\bm{Y}}\rangle_g= F_{\bm{Y}}(f)$ is preserved under $H^1$-convergence. Thus indeed $H^1_{T,\bm{Y}}(\bar{M})$ is an $H^1$-closed subspace and hence a Hilbert space in its own right. We prove the following
\begin{thm}
\label{AppendixAT1}
Let $(\bar{M},g)$ be a compact, connected $3$-manifold and $\bm{Y}\in \mathcal{K}(\bar{M})$. If $H^1_{T,\bm{Y}}(\bar{M})\neq \{0\}$, then there exists some $f\in \left(C^{\infty}(\bar{M})\setminus \{0\}\right)\cap \left(\mathcal{H}^{0}_D(\bar{M}) \right)^{\perp}$ such that
\[
\Delta f:=\delta df=\lambda f\text{, }g(\text{grad}(f),\bm{Y})\equiv 0\text{ and }f|_{\partial\bar{M}}\equiv 0,
\]
for some suitable $\lambda>0$, where $\delta$ denotes the adjoint derivative.
\end{thm}
Our Laplacian is chosen to be a positive operator, i.e. it coincides with the negative standard Euclidean Laplacian in the Euclidean setting. Note further that if $\partial\bar{M}\neq \emptyset$, we have $\left(\mathcal{H}^{0}_D(\bar{M}) \right)^{\perp}=\{0\}^{\perp}=L^2(\bar{M})$.
\newline
\newline
\underline{Proof of \cref{AppendixAT1}:} First of all we equip the space $H^1_{T,\bm{Y}}(\bar{M})$, as defined in (\ref{AppendixA1}) with the norm $\norm{\cdot}_{H^1_0}$ defined by $\norm{f}_{H^1_0}:=\norm{df}_{L^2}$. It follows from \cite[Theorem 2.4.10]{S95} that the norms $\norm{\cdot}_{H^1_0}$ and $\norm{\cdot}_{H^1}$ are equivalent on $H^1_{T,\bm{Y}}(\bar{M})$. We then consider the following constrained minimisation problem
\begin{MPEquation}
\label{AppendixMP1}
\mathcal{E}:\left(H^1_{T,\bm{Y}}(\bar{M}),\norm{\cdot}_{H^1_0}\right)\rightarrow \mathbb{R}, f\mapsto \frac{1}{2}\norm{f}^2_{H^1_0}\rightarrowtail\text{ min }, \frac{1}{2}\norm{f}^2_{L^2}=1.
\end{MPEquation}
We first observe that the condition $H^1_{T,\bm{Y}}(\bar{M})\neq\{0\}$ implies that the set among which we wish to minimise is non-empty, which can be achieved by scaling any fixed non-zero element accordingly. It is now standard to conclude, using the direct method in the calculus of variations and keeping in mind the Sobolev embeddings for compact manifolds, \cite[Theorem 1.3.6]{S95}, that the constrained minimisation problem (\ref{AppendixMP1}) admits a global minimiser. Further one easily checks that the function $\mathcal{E}$ is continuously Fr\'{e}chet-differentiable with respect to the norm $\norm{\cdot}_{H^1_0}$. Similarly the constraint function $\mathcal{L}:H^1_{T,\bm{Y}}(\bar{M})\rightarrow \mathbb{R}$, $f\mapsto\frac{1}{2}\norm{f}^2_{L^2}$ is continuously Fr\'{e}chet differentiable with respect to the same norm and one readily checks that its derivative at a given point $f\in H^1_{T,\bm{Y}}(\bar{M})$ is given by $\mathcal{L}^{\prime}(f)(\phi)=\left(f,\phi \right)_{L^2}$ for any $\phi\in H^1_{T,\bm{Y}}(\bar{M})$. This in particular implies that the map $\mathcal{L}^{\prime}(f):H^1_{T,\bm{Y}}(\bar{M})\rightarrow\mathbb{R}$ is surjective for every $f$ which is not constant zero. It then follows from the Lagrange multiplier rule for Banach spaces, \cite[p.270 Proposition 1]{Z95}, that there exists some $\lambda\in \mathbb{R}$ such that our global minimiser $f$ satisfies
\begin{AppA}
\label{AppendixA2}
\left(df,d\phi\right)_{L^2}=\lambda \left(f,\phi\right)_{L^2}\text{ for all }\phi\in H^1_{T,\bm{Y}}(\bar{M}).
\end{AppA}
Setting $\phi=f$ we obtain $\norm{f}^2_{H^1_0}=\lambda \norm{f}^2_{L^2}=2\lambda$ since $f$ satisfies the specified constraint. Thus $\lambda \geq 0$ and if $\lambda=0$ were true we would find $\norm{f}_{H^1_0}=0$ and consequently $f\equiv 0$, contradicting the constraint which $f$ satisfies. Thus we have $\lambda >0$. Further we observe that $f\in H^1_{T,\bm{Y}}(\bar{M})\subseteq \left(\mathcal{H}^0_D(\bar{M})\right)^{\perp}$. It therefore follows that $f$ admits a Dirichlet potential $f_D\in H^3(\bar{M})$, \cite[Theorem 2.2.4, Theorem 2.2.6]{S95}, which satisfies $f_D|_{\partial\bar{M}}\equiv 0$, $f_D\in \left(\mathcal{H}^0_D(\bar{M})\right)^{\perp}$ and solves the equation
\begin{AppA}
\label{AppendixA3}
\Delta f_D=f.
\end{AppA}
Note that if we can show that $F_{\bm{Y}}(f_D)=0$, then we will have $f_D\in H^1_{T,\bm{Y}}(\bar{M})$. To this end we observe first that for every $f\in C^{\infty}(\bar{M})$ we have the estimate $\norm{F_{\bm{Y}}(f)}_{H^2}\leq C\norm{f}_{H^3}$ for some constant $C>0$ independent of $f$, see \cite[Chapter 1.3]{S95} for an introduction into the theory of higher order Sobolev spaces. We conclude that there is a unique continuous extension of $F_{\bm{Y}}$ to the domain $H^3(\bar{M})$ with range $H^2(\bar{M})$. By uniqueness of the extensions and since $H^3(\bar{M})\subseteq H^1(\bar{M})$ the linear extensions to the domains $H^1(\bar{M})$ and $H^3(\bar{M})$ must coincide on $H^3(\bar{M})$. Since $f_D\in H^3(\bar{M})$ this implies $F_{\bm{Y}}(f_D)\in H^2(\bar{M})$. Further by definition of the space $H^3(\bar{M})$ there exists a sequence $(f_n)_n\subset C^{\infty}(\bar{M})$ which converges to $f_D$ in $H^3$ and consequently $\Delta f_n$ converges to $\Delta f_D$ in $H^1$. Thus on the one hand $F_{\bm{Y}}(\Delta f_n)$ converges to $F_{\bm{Y}}(\Delta f_D)$ in $L^2$ and on the other hand, $\Delta F_{\bm{Y}}(f_n)$ converges in $L^2$ to $\Delta F_{\bm{Y}}(f_D)$ since $F_{\bm{Y}}(f_n)$ converges to $F_{\bm{Y}}(f_D)$ in $H^2$. However, one easily checks by means of the Killing equations that for every $\phi\in C^{\infty}(\bar{M})$ we have $\Delta F_{\bm{Y}}(\phi)=F_{\bm{Y}}(\Delta \phi)$. Thus we obtain, since $F_{\bm{Y}}(f)=0$, and by means of (\ref{AppendixA3})
\[
0=F_{\bm{Y}}(f)=F_{\bm{Y}}(\Delta f_D)=\lim_{n\rightarrow \infty}F_{\bm{Y}}(\Delta f_n)=\lim_{n\rightarrow \infty}\Delta F_{\bm{Y}}(f_n)=\Delta F_{\bm{Y}}(f_D).
\]
It follows further from the Sobolev embeddings, \cite[Theorem 1.3.6]{S95}, that $H^3(\bar{M})$ embeds continuously into $C^1(\bar{M})$, so that in fact $f_D\in C^1(\bar{M})$. Now recall that $f_D|_{\partial\bar{M}}\equiv 0$ and therefore the tangent part of its gradient must vanish, i.e. $\text{grad}(f_D)\perp \partial\bar{M}$. This implies, since $\bm{Y}$ is tangent to the boundary, that $F_{\bm{Y}}(f_D)=g(\bm{Y},\text{grad}(f_D))$ vanishes identically on the boundary. Integration by parts then implies $0=\left(\Delta F_{\bm{Y}}(f_D),F_{\bm{Y}}(f_D) \right)_{L^2}=\left(dF_{\bm{Y}}(f_D),dF_{\bm{Y}}(f_D) \right)_{L^2}$, where the boundary term vanishes due to the fact that $F_{\bm{Y}}(f_D)$ vanishes on the boundary. Hence $dF_{\bm{Y}}(f_D)\equiv 0$. It then follows from the boundary conditions which $F_{\bm{Y}}(f_D)$ satisfies, that we have $F_{\bm{Y}}(f_D)\in \mathcal{H}^0_D(\bar{M})$, i.e. that $F_{\bm{Y}}(f_D)$ is constant. If $\partial\bar{M}\neq \emptyset$, then $\mathcal{H}^0_D(\bar{M})=\{0\}$, which implies that $F_{\bm{Y}}(f_D)\equiv 0$. If instead $\partial\bar{M}=\emptyset$, then the $C^1$ function $f_D$ must by compactness attain a global maximum on $\bar{M}$ and since all points are interior points and due to the $C^1$-regularity, we see that $\text{grad}(f_D)$ vanishes at at least one point. But this implies that $F_{\bm{Y}}(f_D)=g(\bm{Y},\text{grad}(f_D))$ has a zero. Since $F_{\bm{Y}}(f_D)$ is constant, it must be identically zero in this case as well. We conclude that in any case we have $F_{\bm{Y}}(f_D)=0$ and overall $f_D\in H^1_{T,\bm{Y}}(\bar{M})$. Now since $f_D$ is the Dirichlet potential of $f$ it in particular satisfies
\[
\left(df_D,d\eta\right)_{L^2}=\left(f,\eta\right)_{L^2}\text{ for very }\eta\in H^1(\bar{M})\text{ with }\eta|_{\partial\bar{M}}=0.
\]
Since all elements of $H^1_{T,\bm{Y}}(\bar{M})$ obey Dirichlet boundary conditions, the above equality in combination with (\ref{AppendixA2}) implies
\[
\left(d\left(\frac{f}{\lambda}-f_D\right),d\phi \right)_{L^2}=0\text{ for all }\phi\in H^1_{T,\bm{Y}}(\bar{M}).
\]
We have shown that $f_D\in H^1_{T,\bm{Y}}(\bar{M})$ and thus we may choose $\phi=\frac{f}{\lambda}-f_D$ to conclude that $\norm{\frac{f}{\lambda}-f_D}_{H^1_0}=0$, which in turn implies $f=\lambda f_D$. Now a standard bootstrapping argument by means of Sobolev embeddings, \cite[Theorem 1.3.6]{S95}, and the regularity properties of the Dirichlet potential, \cite[Theorem 2.2.6]{S95}, imply that $f$ is smooth. Further (\ref{AppendixA3}) yields $\Delta f=\lambda f$, where $\lambda>0$. Also $f$ is not identically zero, as it is a global minimiser of the minimisation problem (\ref{AppendixMP1}). This concludes the proof. $\square$
\begin{rem}
\label{AppendixAR2}
If $(\bar{M},g)$ is real analytic, then by standard elliptic regularity results the eigenfunction $f$ is real analytic as well.
\end{rem}
The following ansatz in order to construct Beltrami fields from the solutions of the symmetry constrained Laplacian Dirichlet eigenproblem is inspired by a preprint of Gavrilov\footnote{arXiv identifier: 1906.07465} [Section 5.3] and \cite[Chapter 4, Proposition 3,4]{C99}, where the rotationally symmetric Euclidean case was considered.
\begin{cor}
\label{AppendixAC3}
Let $(\bar{M},g)$ be a compact, connected $3$-manifold and suppose that $\bm{Y}\in \mathcal{K}(\bar{M})\setminus\{0\}$ satisfies $\text{curl}(\bm{Y})=\kappa \bm{Y}$ for some $\kappa\in \mathbb{R}$. If $H^1_{T,\bm{Y}}(\bar{M})\neq \{0\}$, then there exists some $f\in C^{\infty}(\bar{M})$ obeying $\bm{Y}(f)\equiv 0$, $f|_{\partial\bar{M}}\equiv 0$ and some $\mu>0$ such that the vector field $\bm{X}:=\bm{Y}\times \text{grad}(f)-f\mu\bm{Y}$ satisfies
\begin{AppA}
\text{curl}(\bm{X})=\mu\bm{X},\text{ }[\bm{Y},\bm{X}]\equiv 0\text{ and }\bm{Y}\times \bm{X}\not\equiv 0.
\end{AppA}
In particular $\bm{X}$ is tangent to the boundary. Further $f$ can be chosen real analytic, whenever $(\bar{M},g)$ is real analytic.
\end{cor}
\underline{Proof of \cref{AppendixAC3}:} Since $H^1_{T,\bm{Y}}(\bar{M})\neq \{0\}$, we obtain from \cref{AppendixAT1} some $\lambda>0$ and a smooth, non-constant function $f$ obeying Dirichlet boundary conditions and satisfying $\Delta f=\lambda f$ and $\bm{Y}(f)\equiv 0$. Define $\mu:=\frac{\kappa}{2}+\sqrt{\frac{\kappa^2}{4}+\lambda}>0$, which is well defined since $\lambda >0$. In particular $\mu$ solves the quadratic equation $\mu^2-\mu \kappa-\lambda=0$. It follows from \cref{MP7} that $g(\bm{Y},\bm{Y})$ is constant and since we assume $\bm{Y}$ not to be the zero vector field, we have $g(\bm{Y},\bm{Y})\equiv c>0$. Thus upon rescaling $\bm{Y}$ by some suitable constant, we may assume that $c=1$. Now define $\bm{X}$ as proposed in the statement of the corollary and compute
\[
\text{curl}(\bm{X})=\text{curl}(\bm{Y}\times \text{grad}(f))-\mu \left(\text{grad}(f)\times \bm{Y}\right)-\mu f\kappa \bm{Y},
\]
where we used that $\bm{Y}$ is a Beltrami field. By standard calculus identities and since $\bm{Y}$ is divergence-free we obtain $\text{curl}(\bm{Y}\times \text{grad}(f))=-(\Delta f)\bm{Y}-[\bm{Y},\text{grad}(f)]$, where we recall that $\Delta$ is chosen to be a positive operator, hence the minus sign. By the Killing properties of $\bm{Y}$, \cref{PT1L1}, we find $[\bm{Y},\text{grad}(f)]=\text{grad}(g(\bm{Y},\text{grad}(f)))\equiv 0$, since $g(\bm{Y},\text{grad}(f))=\bm{Y}(f)\equiv 0$ by choice of $f$. Since $\Delta f=\lambda f$ we arrive at
\[
\text{curl}(\bm{X})=\mu\left(\bm{Y}\times \text{grad}(f)-\bm{Y}f\left(\kappa+\frac{\lambda}{\mu} \right) \right)=\mu \left(\bm{Y}\times \text{grad}(f)-f\mu \bm{Y} \right)=\mu\bm{X},
\]
by definition of $\mu$. This in particular implies that $\bm{X}$ is divergence-free and consequently $[\bm{Y},\bm{X}]=\text{curl}(\bm{X}\times \bm{Y})$. By standard calculus identities and definition of $\bm{X}$ we find
\[
\bm{X}\times \bm{Y}=\left(\bm{Y}\times\text{grad}(f)\right)\times \bm{Y}=\text{grad}(f)g(\bm{Y},\bm{Y})-\bm{Y}g(\bm{Y},\text{grad}(f))=\text{grad}(f),
\]
where we used that $g(\bm{Y},\text{grad}(f))=\bm{Y}(f)\equiv 0$ by properties of $f$ and that $g(\bm{Y},\bm{Y})\equiv 1$, by properties of $\bm{Y}$ and our scaling. This immediately implies that $\bm{X}$ and $\bm{Y}$ commute and further it shows that if $\bm{X}\times \bm{Y}\equiv 0$ then $f$ is constant, which contradicts our choice of $f$. The regularity assertion follows from \cref{AppendixAR2}. $\square$
\begin{rem}
Observe that the function $f$ in the construction of $\bm{X}$ in \cref{AppendixAC3} vanishes on the boundary. This implies that the tangent part of its gradient vanishes and hence its gradient is normal to the boundary. Therefore $\bm{X}$ is tangent to the boundary. Further this implies $\bm{X}=\bm{Y}\times \text{grad}(f)$ on $\partial\bar{M}$, which shows that $\bm{X}$ and $\bm{Y}$ are by construction everywhere $g$-orthogonal on $\partial\bar{M}$. This confirms the reasoning of item (iv) of \cref{MR2}.
\end{rem}
The following shows that all commuting Beltrami fields can be obtained from the procedure described in the proof of \cref{AppendixAC3}
\begin{prop}
\label{AppendixAP5}
Let $(M,g)$ be a compact, connected $3$-manifold with empty boundary. Suppose we are given $\bm{Y}\in \mathcal{K}(M)\setminus \{0\}$ satisfying $\text{curl}(\bm{Y})=\kappa \bm{Y}$ for some $\kappa\in \mathbb{R}$. If $\bm{X}\in \mathcal{V}(M)$ obeys
\[
\text{curl}(\bm{X})=\mu\bm{X}\text{ and }[\bm{Y},\bm{X}]\equiv 0\text{ for some }\mu\in \mathbb{R}\setminus\{0\},
\]
then there exists a smooth function $f\in C^{\infty}(M)$ with
\[
\bm{X}=\bm{Y}\times \text{grad}(f)-\mu f\bm{Y},\text{ }\bm{Y}(f)\equiv 0\text{ and }\Delta f=\mu(\mu-\kappa)f.
\]
If in addition $\bm{Y}\times \bm{X}\not\equiv 0$, then $f$ is not constant, in particular $H^1_{T,\bm{Y}}(M)\neq \{0\}$, and $\mu(\mu-\kappa)>0$. 
\end{prop}
\underline{Proof of \cref{AppendixAP5}:} It follows from \cref{MP7} that $g(\bm{Y},\bm{Y})\equiv c>0$ is constant, since $\bm{Y}$ is not the zero vector field. We then define $f:=-\frac{g(\bm{Y},\bm{X})}{c\mu}$ and observe that by means of \cref{PT1L1} and properties of $\bm{X}$, we have
\begin{AppA}
\label{AppendixA4}
-c\mu\text{ }\text{grad}(f)=\text{grad}(g(\bm{Y},\bm{X}))=\mu\bm{Y}\times \bm{X}.
\end{AppA}
We see that $\text{grad}(f)$ is $g$-orthogonal to $\bm{Y}$ and consequently $\bm{Y}(f)=g(\bm{Y},\text{grad}(f))\equiv 0$ as claimed. Further taking the cross product with $\bm{Y}$ from the left of the above identity and using the triple vector product rule, we obtain
\[
-c\left(\bm{Y}\times \text{grad}(f)\right)=\bm{Y}\times \left(\bm{Y}\times \bm{X} \right)=\bm{Y}g(\bm{Y},\bm{X})-g(\bm{Y},\bm{Y})\bm{X}\Leftrightarrow \bm{X}=\bm{Y}\times \text{grad}(f)-\mu f\bm{Y},
\]
where we used the definition of $f$ and $g(\bm{Y},\bm{Y})\equiv c$. Using the identity $\Delta=-\text{div}(\text{grad})$ we compute by (\ref{AppendixA4})
\[
c\Delta f=\text{div}(\bm{Y}\times \bm{X})=g(\text{curl}(\bm{Y}),\bm{X})-g(\text{curl}(\bm{X}),\bm{Y})=(\kappa-\mu)g(\bm{X},\bm{Y})
\]
by properties of $\bm{Y}$ and $\bm{X}$ and by standard calculus identities. The definition of $f$ implies $\Delta f=\mu(\mu-\kappa)f$ as claimed. Now if $f$ is constant we conclude from (\ref{AppendixA4}) that $\bm{Y}\times \bm{X}\equiv 0$. Lastly since $\Delta$ is a positive operator we have $\mu(\mu-\kappa)\geq 0$, while $\mu(\mu-\kappa)=0$, implies $\Delta f=0$, which in turn on compact, boundaryless manifolds implies that $f$ is constant. $\square$
\section{Killing flows in lower dimensions}
\label{Appendix}
In this section we assume only the Hausdorff and second countability property by default, in particular we do not assume orientability. All additional requirements will be explicitly stated. The key to the characterisation in lower dimensions is an analogue of \cref{PT1L1} for the special case $\bm{X}=\bm{Y}$
\begin{lem}
\label{AL1}
Let $(\bar{M},g)$ be a smooth, Riemannian manifold of dimension $n$, with or without boundary, and let $\bm{Y}\in \mathcal{K}(\bar{M})$ be a smooth Killing field. Then
\[
g(\bm{Y},\text{grad}(g(\bm{Y},\bm{Y})))=0\text{ and }\nabla_{\bm{Y}}\bm{Y}=-\frac{1}{2}\text{grad}(g(\bm{Y},\bm{Y})).
\]
\end{lem}
The proof of this lemma follows from direct calculations, keeping in mind the Killing equations. In particular the conclusions are still true if $\bm{Y}$ is not tangent to the boundary, but only satisfies the Killing equations. Now recall that being tangent to the boundary in dimension $1$ means to vanish at the boundary, which allows for the following characterisation
\begin{prop}[Characterisation of smooth Killing flows in dimension $1$]
\label{AP2}
Let $(\bar{M},g)$ be a connected, smooth, Riemannian manifold of dimension $1$, with or without boundary. Then we have
\begin{enumerate}
\item If $\partial\bar{M}\neq \emptyset:$ $\mathcal{K}(\bar{M})=\{0\}$.
\item If $\partial\bar{M}=\emptyset$, then $\bm{Y}\in \mathcal{K}(\bar{M})$ if and only if $\bm{Y}$ is either the zero vector field or, after rescaling by some constant factor $c>0$, is a smooth vector field whose field lines are all unit speed geodesics.
\end{enumerate}
\end{prop}
\underline{Proof of \cref{AP2}:} In dimension $1$ we have $2\nabla_{\bm{Y}}\bm{Y}=\text{grad}(g(\bm{Y},\bm{Y}))$ and so if $\bm{Y}$ is a Killing field, then \cref{AL1} implies $\text{grad}(g(\bm{Y},\bm{Y}))\equiv 0$ and in conclusion $\nabla_{\bm{Y}}\bm{Y}\equiv 0$. This proves the first direction in (ii) and if $\partial\bar{M}$ is non-empty, then the tangent to the boundary condition, implies that $\bm{Y}=0$ at the boundary and hence we must have $g(\bm{Y},\bm{Y})\equiv 0$, which proves (i). As for the converse of (ii), we either have $\bm{Y}\equiv 0$ which obviously is a Killing field, or all field lines of $\bm{Y}$ (up to rescaling by a constant) are unit speed geodesics. This implies $\nabla_{\bm{Y}}\bm{Y}=0$ and that $\bm{Y}$ is no-where vanishing. Writing out $\nabla_{\bm{Y}}\bm{Y}$ in normal coordinates centred around some fixed point $p\in \bar{M}$ yields $0=Y^1(p)(\partial_1Y^1)(p)\partial_1(p)$ and since $Y(p)\neq 0$ we have $(\partial_1Y^1)(p)=0$, which are exactly the Killing equations in dimension $1$ in normal coordinates. Since $p$ was arbitrary we conclude that $\bm{Y}$ is a Killing field. $\square$
\newline
\newline
We now come to the $2$-dimensional case
\begin{prop}[Characterisation of real analytic Killing flows in dimension $2$]
\label{AP3}
Let $(\bar{M},g)$ be a compact, connected, real analytic, Riemannian manifold of dimension $2$, with or without boundary, and let $\bm{Y}\in \mathcal{K}^{\omega}(\bar{M})$. Then one of the following $3$ situations occurs
\begin{enumerate}
\item $\bm{Y}\equiv 0$.
\item After rescaling $\bm{Y}$ by a suitable constant factor $c>0$ all field lines of $\bm{Y}$ are unit speed geodesics.
\item There exists a compact, $\mathcal{H}^1$-countably-$1$ rectifiable subset $\Gamma\subset \bar{M}$ in the sense of Federer, which contains the boundary of $\bar{M}$ and such that $\bar{M}\setminus \Gamma$ is the disjoint union of the images of integral curves of $\bm{Y}$, each of which is a real analytically embedded circle $S^1$ in the interior of $\bar{M}$. Further, the connected component of $\bar{M}\setminus \Gamma$ containing a given invariant circle is an open neighbourhood of the circle and a union of invariant circles. In particular every connected component of $\bar{M}\setminus \Gamma$ is diffeomorphic to $\mathbb{R}\times S^1$.
\end{enumerate}  
\end{prop}
\underline{Proof of \cref{AP3}:} If $f:=g(\bm{Y},\bm{Y})$ is constant, it follows from \cref{AL1} that $\nabla_{\bm{Y}}\bm{Y}\equiv 0$ and hence all field lines of $\bm{Y}$ are constant speed geodesics of the same speed, which corresponds to the first two cases. Now if $f$ is not constant, then we can argue similarly as in the proof of \cref{PT3P1}, as long as we make sure that $\text{grad}(f)$ is normal to the boundary. However since $\bm{Y}$ is tangent to the boundary, we either have $\bm{Y}(p)\neq0$ at some given point $p\in \partial\bar{M}$, which due to the pointwise $g$-orthogonality of $\bm{Y}$ and $\text{grad}(f)$, \cref{AL1}, implies that the gradient of $f$ must be normal to the boundary (since the tangential space is $1$-dimensional) or we have $\bm{Y}(p)=0$, which implies $\nabla_{\bm{Y}}\bm{Y}(p)=0$, which in turn by \cref{AL1} implies that the gradient of $f$ vanishes at $p$ and hence is normal to the boundary. We can now define $\Gamma$ in correspondence with the proof of \cref{PT3P1} and just like in the proof of the latter, the set $\bar{M}\setminus \Gamma$ decomposes into the connected components of the compact, regular level sets of $f|_M$, with $M=\text{int}(\bar{M})$, which are all compact, $1$-dimensional, real analytically embedded submanifolds without boundary and hence must be diffeomorphic to $S^1$. Observe that the gradient of $f$ does not vanish on $\bar{M}\setminus \Gamma$ and so $\bm{Y}$ does not either, because otherwise this would imply that $\nabla_{\bm{Y}}\bm{Y}$ vanishes, contradicting the fact that the gradient of $f=g(\bm{Y},\bm{Y})$ does not vanish. Since $\bm{Y}$ is $g$-orthogonal to $\text{grad}(f)$ it must be tangent to the level sets of $f$. Now the field lines of $\bm{Y}$ are non-constant and travel with some minimal, strictly positive speed along the connected components of these level sets. Hence the image of each field line starting at such an invariant circle must be the whole circle. The statement about the connected components follows exactly as in the proof of \cref{PT3P1}. $\square$
\section*{Acknowledgements}
This work has been funded by the Deutsche Forschungsgemeinschaft (DFG, German Research Foundation) – Projektnummer 320021702/GRK2326 –  Energy, Entropy, and Dissipative Dynamics (EDDy). I would like to thank Christof Melcher and Heiko von der Mosel for discussions. Further I want to thank Daniel Peralta-Salas for pointing out Gavrilov's work to me.
\bibliographystyle{plain}
\bibliography{mybibfile}
\footnotesize
\end{document}